\documentclass[12pt,a4paper]{article}

\usepackage
{hyperref} %colora i riferimenti

\usepackage{amsmath,amsfonts,amsthm}

\usepackage{graphicx}
\usepackage{latexsym}
\usepackage{enumitem}
\usepackage{xcolor}

\usepackage{caption}
\usepackage[utf8]{inputenc}
\usepackage{orcidlink}

\usepackage[title,titletoc]{appendix}

\newcommand{\Rset}{\mathbb R}
\newcommand{\D}{{\bf D}_0}

\newtheorem{theorem}{Theorem}[section]
\newtheorem{lemma}[theorem]{Lemma}
\newtheorem{corollary}[theorem]{Corollary}
\newtheorem{definition}[theorem]{Definition}

\newtheorem{remark}[theorem]{Remark}
\newtheorem{proposition}[theorem]{Proposition}
\newtheorem*{theorem*}{Theorem}

\newcommand{\abs}[1]{\left |#1\right |}
\newcommand{\abss}[1]{\left |#1\right |_*}

\newcommand{\xg}{a_1}
\newcommand{\xp}{a_2}
\newcommand{\decl}{:=}
\newcommand{\HS}{\Rset^N_+}
\newcommand{\DHS}{\partial\Rset^N_+}
\newcommand{\CHS}{\overline{\HS}}
\newcommand{\ringentire}{(${\cal R}$)\ }
\newcommand{\rhh}{$({\cal R}^+)$ }
\newcommand{\rhz}{$({\cal R}^+_0)$ }
\newcommand{\C}{{\cal C}}

\begin{document}

\title{Characterization of positive superharmonic\\
functions in a half-space 
%\thanks{We are very pleased to acknowledge the support of FRA 2024 - Universit\`a degli Studi di Trieste.} 
}

\author{Lorenzo D'Ambrosio\thanks{Universit\`{a}  degli Studi  di Udine, 
Dipartimento di Scienze Matematiche, Informatiche e Fisiche,
Via delle Scienze 206, I-33100 Udine, Italy. Member of INdAM-GNAMPA. \orcidlink{0000-0003-0677-056X} 
email: \href{mailto:lorenzo.dambrosio@uniud.it}{\textcolor{black}{\tt{lorenzo.dambrosio@uniud.it}}} } \ and Enzo Mitidieri\thanks{Universit\`a degli Studi di Trieste,
Dipartimento di Matematica, Informatica e Geoscienze,
via A.Valerio, 12/1, I-34127 Trieste, Italy, \orcidlink{0000-0001-5042-9401}
email: \href{mailto:mitidier@dmi.units.it}{\textcolor{black}{\tt mitidier@dmi.units.it  }}}}

\renewcommand{\theenumi}{\Alph{enumi}}

\maketitle

{\abstract{We prove a representation formula for superharmonic functions on the 
half-space $\Rset^N_+:=\Rset^{N-1}\times]0,+\infty[$.
As a consequence, we derive some comparison principles and various positivity results.}

\medskip


\medskip

\noindent{\bf AMS MSC(2020):} 
%35C99, %(1973-now) Nonlinear higher-order PDEs 
35C15, %Integral representations of solutions to PDEs
%28C05, % (1973-now) Biharmonic and polyharmonic equations and functions in higher dimensions 
%35B53, % (2010-now) Liouville theorems and Phragmén-Lindelöf theorems in context of PDEs 
%35J61.% Semilinear elliptic equations
35R06, %   	PDEs with measure
35B45. %  	A priori estimates in context of PDEs

\tableofcontents
\section{Introduction}
 Let $u\in L^1_{loc}(\Rset^N)$, $N\ge 3$, be a nonnegative superharmonic function,
namely,  
\begin{equation}\label{eq:entire}
-\Delta u=\mu,\quad u\ge 0,\qquad  in\ \  {\cal D}'(\Rset^N),
\end{equation}
here $\mu$ is a positive Radon measure.\footnote{As usual, by a positive
 \emph{Radon measure}, we mean a regular Borel measure which is positive and
 finite on compact sets.}
From Riesz representation theorem (see \cite[Corollary 4.4.2]{AG}) it follows that $u$ is given by,
\begin{equation}\label{repEntire}
  u(x)=l+\int_{\Rset^N} \Gamma^x(y) d\mu(y), \quad for\ a.e. \ x\in\Rset^N,
\end{equation}
where $\Gamma^x(y)$ 
is the fundamental solution of $-\Delta$ (see \eqref{def:greenD}) and $l=\inf u$.

The assumption of nonnegativity of u is crucial for the validity of \eqref{repEntire}. Therefore, finding sufficient conditions that ensure the positivity of a superharmonic function is an important and natural question.
An answer to this problem is provided in \cite{CDM}, where the authors demonstrate:
\begin{theorem} \label{th:RepClassic}
Let $N\ge 3$ and let $u\in L^1_{loc}(\Rset^N)$ be a superharmonic function.
Let $\mu=-\Delta u$ and $l\in \Rset$.
Then  \eqref{repEntire} holds if and only if
\begin{equation*}
\liminf_{R\to+\infty}\frac{1}{R^N}\int_{R<|x-y|<2R} |u(y)-l|\,dy=0\quad for\ a.e. \ x\in \Rset^N.\eqno{({\cal R})}
\end{equation*}
\end{theorem}

The ring condition \ringentire also plays a role in the representation theorem related to distributional solutions of $L(u)\ge 0$, even when $L$ is a higher-order operator, either homogeneous or non-homogeneous (see \cite{CDM}, \cite{DG22}). Representation theory has also been studied for problems stated in non-Euclidean settings (see \cite{BLU}, \cite{DM19}, and \cite{DMP}). Further generalizations to smooth solutions of problems associated with general second-order operators $L(u)\ge 0,$ have been studied in \cite{DM19}, where the role of \ringentire is played by special means modeled on the fundamental solution of the operator $L(u)\ge 0$.

One  purpose of this paper  is to study superharmonic functions on the half-space 
$\HS:=\{(x_1,\dots,x_N)\in\Rset^{N}, x_N>0 \}$ where $N\ge 2,$  and related representation formulae.  Additionally, we will prove some comparison principles (see \cite{damiII} for some Liouville theorems). 
Here, we will focus on  problems in  half-space $\HS;$ the more general case of a cone in $\Rset^N$ will be addressed  in a forthcoming paper.

It is well known that  for nonnegative superharmonic functions in the half-space $\HS$, the following Riesz representation theorem holds. See Corollary 4.4.2 and Theorem 4.3.8 in \cite{AG}.
\begin{theorem}[Riesz representation]\label{th:RepHSClassic}
 Let $u$ be a
distributional superharmonic function in $\HS$, that is $-\Delta u\ge 0$ in ${\cal D}'(\HS)$. 
  If $u\ge 0$, then there exist a positive Radon measure
  $\mu$ on $\HS$ and a positive Radon measure $\nu$ on $\DHS$ together with a constant $h\ge 0,$ such that
  \begin{equation}\label{RepHSClassic} 
    u(x)=h x_N+ \int _{\DHS} K^x(y') d\nu(y') +  \int_{\HS} G^x(y) d\mu(y)\quad for\ a.e. \ x\in \HS.
  \end{equation}
  Furthermore
  $$-\Delta u=\mu \ \ on\ \HS\ \ in\  the\  distributional\ sense.$$
\end{theorem}
Here $K$ and $G$ are the Poisson kernel and the Green function of $-\Delta $ on $\HS$ respectively 
(see Appendix \ref{sec2} for more details).

Even in the context of problems in a half-space,  the question of the nonnegativity of superharmonic functions arises.

Throughout this paper, unless otherwise specified,  we will assume that $\mu$ and $\nu$ are
 positive Radon measures on $\Rset^N_+$ and on  $\DHS=\Rset^{N-1}$ respectively.
Our research will focus on the differential equation, 
\begin{equation}\label{eqh}
 \left\{\begin{array}{ll}
-\Delta u = \mu& \mathrm{on }\ \HS, \cr \cr
u= \nu & \mathrm{on }\ \DHS, \cr
\end{array}
 \right.
\end{equation}
and on the differential inequality
\begin{equation}\label{dish}
 \left\{\begin{array}{ll}
-\Delta u \ge \mu & \mathrm{on }\ \HS, \cr \\
u \ge \nu & \mathrm{on }\ \DHS. \cr
\end{array}
 \right.
\end{equation}

As will become increasing evident as we proceed, when studying problems in $\HS$, the role of the ring condition \ringentire for the problem \eqref{eq:entire} in the whole space $\mathbb{R}^N$, in $\HS$ it is played by  the following weighted ring condition:

\begin{equation*}
\begin{aligned}
 there\ exists\ h\in\Rset\ such\ that\qquad\qquad\qquad\qquad \\ 
 \liminf_{R\rightarrow +\infty} \frac{1}{R^{N+2}}\int_{\{y_N>0\} \cap \{R<\abss{x-y}<2R\}} y_N |u(y)-h y_N| \ dy =0,\  for\ a.e.\ x\in\Rset^N_+,
\end{aligned} \eqno{({\cal R}^+)}
\end{equation*}
where $\,\abss\cdot\, $ is the cylindrical norm defined by,
\begin{equation}\label{def:norm}
|x|_\ast=|(x',x_N)|_\ast\decl \max \{|x'|,|x_N|\}.\end{equation}

To understand why, we introduce the condition \rhh\!\!, we first note that  for the problem \eqref{eq:entire}, the solution $u$ is not unique. Indeed, adding a constant  to $u,$ yields a different solution. This non-uniqueness also applies to problem \eqref{eqh}, where a linear
term like $hx_N,\, h\in\Rset,$ can be added to a solution to produce another solution. 
This explains why a condition like \rhh\!\! is necessary for problem \eqref{eqh}.  

 Our first main result is the following.
\begin{theorem*}[\bf A]\label{teo:repintro}Let $N\geq 2$, let $\mu$ and $\nu$ be 
 positive Radon measures on $\Rset^N_+$ and  $\DHS,$ respectively. Assume that
 $\mu\in {\cal M}_{loc}(\CHS,x_N)$ and let $u\in L^1_{loc}(\CHS)$  be a weak solution of (\ref{dish}).
Suppose $x\in\HS$ is a Lebesgue point for $u$. Assume that there exists $h\in\Rset$ such that  \rhh holds at the point $x$ that is,
\begin{equation*}
 \liminf_{R\rightarrow +\infty} \frac{1}{R^{N+2}}\int_{\{y_N>0\} \cap \{R<\abss{x-y}<2R\}} y_N |u(y)-h y_N| \ dy =0.
\end{equation*}
   Then,
   \begin{equation}\label{intro:disrepmeas+}
      u(x)\ge hx_N+\int _{\DHS} K^x(y') d\nu(y') +  \int_{\HS} G^x(y) d\mu(y).
   \end{equation}
   
   Moreover, if $u$ is a weak solution of \eqref{eqh}, then the inequality in \eqref{intro:disrepmeas+} becomes  an equality, that is
     \begin{equation}\label{intro:eqrepmeas+}
      u(x)= hx_N+\int _{\DHS} K^x(y') d\nu(y') +  \int_{\HS} G^x(y) d\mu(y).
   \end{equation}
\end{theorem*}

We observe that a special case arises  when \rhh holds with $h=0$, i.e., 
\begin{equation*}
 \liminf_{R\rightarrow +\infty} \frac{1}{R^{N+2}}
  \int_{\{y_N>0\} \cap \{R<\abss{x-y}<2R\}} y_N |u(y)| \ dy =0\  for\ a.e.\ x\in\Rset^N_+.\eqno{({\cal R}^+_0)}
\end{equation*}

Notice that all bounded and measurable functions on $\HS$, satisfy \rhz\!\!.  
Other examples of functions satisfying \rhz\!\!, will be provided in Appendix~\ref{sec:R0} below.

We point out that another interesting class of functions satisfying \rhz is the set of solutions to the problem,
\begin{equation}
 \left\{\begin{array}{ll}
\pm \Delta u \ge |u|^q, & \mathrm{on }\ \HS, \cr
u = 0, & \mathrm{on }\ \DHS, \cr
\end{array}
 \right.
\end{equation}
here $q>1.$ For a proof of this fact and related Liouville theorems see \cite{damiII}.

Our main 
result Theorem (A), can also be applied to the representation of harmonic functions.
Recently, in \cite{zha16},
integral representation formulae for harmonic functions on a half-space have been studied.
It is worth comparing \rhz\!\!, with the condition given in \cite{zha16}
for representing a harmonic function by an integral formula. 

\medskip
We emphasize that  for a superharmonic function $u$ in a half-space $\HS,$ the condition
\rhh  is sufficient  for the validity of the representation formula \eqref{RepHSClassic}, thereby ensuring its  positivity, provided $h\ge0$.
Furthermore, the condition \rhh is also necessary for the validity of the integral representation. Indeed, our second main result is the following theorem whose principal statement is the converse of the statement in Theorem (A).

\begin{theorem*}[\bf B] Let $N\geq 2$,  and let $\mu$ and $\nu$ be 
 positive Radon measures on $\Rset^N_+$ and  $\DHS,$ respectively.
  Let $u$ be defined by \eqref{intro:eqrepmeas+} and assume that it is finite  for a.e. $x\in\HS$.

  Then $u\in L^1_{loc}(\CHS)$, $\mu\in {\cal M}_{loc}(\CHS,x_N)$ and $u$ is a weak solution of \eqref{eqh}.  Moreover, $u$ satisfies
  the condition \rhh\!, where the constant $h$ appears in \eqref{intro:eqrepmeas+}.
  
  Furthermore, 
  \begin{enumerate}
     \item[(a)] $$h=\inf_{\HS} \frac{u(x)}{x_N}$$
  \item[(b)] Let $\Omega\subset\HS$ be a nonempty bounded open set. If $h=\inf_{\Omega} \frac{u(x)}{x_N}$,
  then $u(x)=h x_N$, with $\mu\equiv 0$ and $\nu\equiv 0$. 
  \item[(c)] If $u-hx_N$ is not identically zero, then there exists a constant $c_0>0$ such that
  \begin{equation}\label{intro:est:ubelow}
  u(x)\ge hx_N+ c_0 \frac{x_N}{1+|x|^N},\quad for \ \ all\ \ x\in\HS.
  \end{equation}
  \item[(d)] For almost every $x\in\HS,$
  \begin{equation*}
 \lim_{R\rightarrow +\infty} \frac{1}{R^{N+2}}\int_{\{y_N>0\} \cap B^*_R(x)} y_N |u(y)-h y_N| \ dy =0.
\end{equation*}
  \end{enumerate}
\end{theorem*}

This paper is organized as follows.

Section \ref{sec3} is dedicated to the defining  weak solutions for problems \eqref{eqh} and \eqref{dish},
even when $\nu$ and $\mu$ are general local Radon measures,
and to the concept of {\it lim-trace}. See Definition~\ref{def:trace}.  To illustrate the implications of this concept, we explore its application in  Section \ref{sec3.1}, where we establish the relationship between  distributional and weak solutions of \eqref{eqh} and \eqref{dish}, and certain specific properties 
(see Remark \ref{rem:counterexample}).
 
 Section \ref{sec:main} formulates and proves the main result of  this paper, namely the representation Theorem \ref{teo:repiff}.
The techniques we use to prove this  representation result for (\ref{eqh}) allow us to deal with measures that are general 
local Radon measures, without requiring positivity. 
See Theorem \ref{repmeasnosign}. 
We investigate the case of superharmonic but not necessarily nonnegative functions $u$ in Theorem  \ref{teo:distr=weak} and Corollary  \ref{cor:linkDW}. Consequently, weak solutions of \eqref{dish} that are bounded from below are nonnegative. See Theorem \ref{teo:bb=0}. This section also includes a simple characterization of a comparison principle (see Theorem \ref{teo:potentialink}). The remainder of the section is devoted to proving other related results.

Lastly, for completeness, in Section  \ref{sec:ringreen}, we highlight another representation theorem for a superharmonic function on $\HS$. The results contained in this section are based on general theorems established in \cite{DM19}.

As a final remark, we would like to point out that if $u$ is a superharmonic function and there exists a superharmonic function $g$ such that
$u\ge -g$, then the results contained in Section \ref{sec:main}, apply to $u+g$.
We emphasize that in this case, by Riesz's decomposition theorem, we know that
there exists an harmonic function $H$ on $\HS$, such that
$$u(x)=H(x)+ \int_{\HS} G^x(y) d\mu(y),$$
where $\mu\decl-\Delta u.$
The harmonic function $H$ can be obtained as a weighted limit of $u$ on some rings.
This construction is detailed in Section \ref{sec:ringreen}, see Theorem \ref{teo:repold}.

In Appendix \ref{sec2}, we present  basic  and important information on
the fundamental solution of $-\Delta$ in $\HS$ and its Poisson kernel and
the related Green function. 
This includes a few precise estimates on related potentials. Appendix \ref{sec:R0} contains
basic examples of functions fulfilling the condition \rhz\!\!.

\medskip

\subsubsection*{Notation}

\begin{itemize}[align=left,leftmargin=0.42in,
parsep=\medskipamount, 
labelsep=0.08in]

\item[$\HS$] The half-space  $\Rset^{N-1}\times]0,+\infty[$, $x=(x',x_N)\in\HS$  with 
$x_N>0$ and $x'\in\Rset^{N-1}$.
\item[$\DHS=\Rset^{N-1}\times \{0\}$] The boundary of $\HS$ that we can identify with $\Rset^{N-1}$.

  \item[ $|\cdot|$] The Euclidean norm.
  \item[$B'_R$] The Euclidean ball in $\Rset^{N-1}$ centered at the origin and radius $R$.
\item[$B^{e,N+2}_R(\xi_0)$] The Euclidean ball in $\Rset^{N+2}$ of radius $R$
centered ad $\xi_0\in\Rset^{N+2}$.
  \item[$|x|_\ast$] The cylindrical norm defined as 
	  $|x|_\ast=|(x',x_N)|_\ast\decl \max \{|x'|,|x_N|\}$.
  \item[$B^*_R(x)$] The cylindrical balls in $\HS$ defined as  $\left\{y\in\HS: \max\{|y_N-x_N|,|y'-x'|\}<R   \right\}$.
%  \item[$A^*_R(x)$] The annulus $A^*_R(x)\decl B^*_{2R}(x)\setminus B^*_R(x).$
  \item[$\C_c(E)$, $\C^2_c(E)$]  The spaces of continuous and $\C^2$ functions respectively whose support is compact and contained in $E$.
  \item[${\D}(\Omega)$]{ Test functions space} used along the paper defined as 
   $${\D}(\Omega)\decl\left\{\varphi_{|\Omega}: \varphi\in \C_c^2(\Rset^N),  
 \varphi=0\ on\ \partial\Omega\right\}.$$  
 \item[${\cal D}'(\HS)$] The space of the distributions on $\HS$.
  \item[$\sigma_N$] {The measure of the unit sphere in $\Rset^N$,}\  $\sigma_N\decl \frac{2\pi^{N/2}}{\Gamma(N/2)} $.
  \item[$C_N, C'_N$] The constants defined as $C_N^{-1}\decl \sigma_N\max\{N-2,1\}$, and $C'_N\decl 2/\sigma_N$.
  \item[$\Gamma^x$] The fundamental solution  of $-\Delta$ at $x$, see \eqref{def:greenD}.
  \item[$G^x$] The Green function of $-\Delta$ at $x$ on $\HS$, see \eqref{def:G}.
  \item[$K^x$] The Poisson kernel, see \eqref{def:K}.
  \item[${\cal M}_{loc}(X)$] {The space of local  Radon measure} on the 
   $\sigma$-compact and locally compact Hausdorff space $X$, that is the the space
    of linear continuous functionals on $\C_c(X)$. 
   Moreover if $\mu\in {\cal M}_{loc}(X)$, then $\mu =  \mu^+ - \mu^-$, with
   $\mu^+$ and $\mu^-$ positive Radon measure and set 
   $\abs\mu \decl \mu^+ + \mu^-$.

  \item[${\cal M}_{loc}(\CHS,x_N)$] The space of local Radon measure 
    $\mu\in {\cal M}_{loc}(\HS)$  such that
 $\int_\Omega x_N d|\mu| <\infty$ for any bounded $\Omega\subset\HS$.

  \item[$L^1_{loc}(\CHS)$]  The space of locally integrable function up to the boundary, that is the space of function $u\in L^1_{loc}(\HS)$ such that
 $u\in   L^1(\Omega)$ for any open bounded set $\Omega\subset\HS$.

  \item[$W^{1,p}_{loc}(\CHS)$] The space of functions $u\in  W^{1,p}_{loc}(\HS)$ 
  such that $u\in W^{1,p}(\Omega)$ for any open bounded set $\Omega\subset\HS$.  
  
  \item[$BV_{loc}(\CHS)$] The space of function $u\in BV_{loc}(\HS)$
    such that $u\in   BV(\Omega)$ for any open bounded set $\Omega\subset\HS$.

 \item[$Tr(u)$]  The lim-trace of $u$, see Definition \ref{def:trace}.

\item[$c$] We will use the symbol $c$ to denote a generic positive constant,
	the value of which may change from line to line, and is not dependent on the solutions of the problem under study and may depend on certain  non-essential parameters.
\end{itemize}

\section{Some qualitative properties of weak and distributional solutions}\label{sec3}
In this section, we present the formal definition of weak solutions for the problem \eqref{eqh}  and
discuss some of its properties.

The definition of weak solution of \eqref{eqh} can be formulated even if $\mu$ and $\nu$ are  more general than positive Radon measures, extending to local local Radon measures.
For a $\sigma$-compact and locally compact Hausdorff space $X,$ a 
local Radon measure $\mu$ is a linear continuous functional on $\C_c(X)$.
The space of all the local Radon measure on $X$ will be denoted by
${\cal M}_{loc}(X)$.

Clearly, ${\cal M}_{loc}(X)$ contains any positive Radon measure as well as any signed Radon measure on $X$. 
Furthermore,  if $\nu$ and $\lambda$ are two positive Radon measures on $X$, then
$\nu-\lambda$ defines a linear functional on $\C_c(X)$, that is, an element of 
 ${\cal M}_{loc}(X)$ (note that in general $\nu-\lambda$ is not a signed Radon measure).
 On the other hand, for any $\mu\in {\cal M}_{loc}(X)$,
there exist two positive Radon measures, $\nu$ and $\lambda,$ such that 
$\mu=\nu-\lambda$. 
As usual, we rewrite $\mu$ in terms of the positive measures
$\mu^+$ and $\mu^-,$ as in the classical Jordan decomposition $\mu=\mu^+-\mu^-$.  Therefore, for the functional $\mu\in {\cal M}_{loc}(X)$ we have the following: 
$$\mu(\phi)=\int_X \phi d\mu^+-\int_X \phi d\mu^-
=:\int_X \phi d\mu\qquad\forall \phi\in \C_c(X).$$
As usual, we set
$\abs\mu \decl \mu^+ + \mu^-$, which defines a classical positive Radon measure.
The interested reader may refer to the relationships between classical measures, local Radon measures and their representation in \cite{bou} and \cite{evagar}.

\medskip
%{\bf Relevant Functional Spaces} \bigskip
In what follows we will need the following spaces:

\bigskip
${\cal M}_{loc}(\CHS,x_N)\decl \{ \mu: \mu$ is a local Radon measure   on $\HS$ such that

\hspace{5cm} $\int_\Omega x_N d|\mu| <\infty$ for any bounded $\Omega\subset\HS$\},

\medskip
and
$${\D}={\D}(\HS)\decl\left\{\phi_{|\HS}: \phi\in \C_c^2(\Rset^N),  \varphi(x',0)=0\right\}.$$

\medskip
%{\bf Definition of Weak and Distributional Solutions}
Next we notice that if, instead of $\mu$ and $\nu$ being measures, they are continuous functions and $u$ is a smooth
solution of \eqref{eqh}, then by multiplying the differential equation in \eqref{eqh} by  $\varphi\in\D $ and integrating by parts, we obtain
the identity:
\begin{equation}\label{def:solpre}
  \int_{\HS}\mu(x)\varphi(x) dx+ \int_{\DHS} \nu(x')\varphi_N(x',0) dx'= 
	\int_{\HS}u (-\Delta \varphi).
\end{equation}
Here, and in the following, we denote by $\varphi_N$ the partial derivative of
the function $\varphi$ with respect to $x_N$, i.e.
$\varphi_N(x)=\partial_N\varphi(x)=\frac{\partial}{\partial x_N}\varphi(x)$.

Similarly, if $u$ solves \eqref{dish},
then by multiplying the differential inequality in \eqref{dish} by a nonnegative $\varphi\in\D, $ we obtain
\begin{equation}\label{def:soldispre}
  \int_{\HS}\mu(x)\varphi(x) dx+ \int_{\DHS} \nu(x')\varphi_N(x',0) dx'\le
	\int_{\HS}u (-\Delta \varphi).
\end{equation}

With these preliminaries, we can now provide the formal definition of weak solutions and discuss their properties.
The relations \eqref{def:solpre} and \eqref{def:soldispre} justify the following definition of weak solution.
\begin{definition} \label{def:genweaksol}
 Let $\mu\in {\cal M}_{loc}(\CHS,x_N)$ and $\nu\in {\cal M}_{loc}(\DHS)$.
A function $u$ is a \em{weak solution of \eqref{eqh} [resp. \eqref{dish}]} if $u\in L^1_{loc}(\CHS)$ and
for every test function $\varphi\in \D$ [resp.  $\varphi\in \D$ and $\varphi\ge 0$], there holds,
\begin{equation}\label{eq:defsol}
  \int_{\HS}\varphi(x) d\mu(x)+ \int_{\DHS} \varphi_N(x',0) d\nu(x')= [\le ]  
	\int_{\HS}u (-\Delta \varphi).
\end{equation}
\end{definition}

\begin{remark}\label{rem:2}
Observe that if, in the definition of weak solution, we replace $\mu$ with a general local Radon measure,
  $\mu\in {\cal M}_{loc}(\HS)$, then it is not guaranteed that the terms in
  \eqref{eq:defsol} are all finite.
On the other hand, if  $\mu\in {\cal M}_{loc}(\HS)$ is a positive Radon measure, it is straightforward to recognize that requiring the term
$ \int_{\HS}\varphi(x) d\mu(x)$ to be finite 
for any $\varphi\in\D$ is equivalent to the fact that the measure $\mu$ satisfies
$\int_{\Omega}x_N d\mu(x)<\infty$ on every bounded open set
$\Omega\subset\HS$, i.e., that $\mu\in {\cal M}_{loc}(\CHS,x_N)$.
\end{remark}
\begin{remark}\label{rem:weak>distr}
Needless to say, if $u$ is a weak solution of \eqref{eqh} then 
$$-\Delta u =\mu \quad in\quad {\cal D}'(\HS),$$ 
i.e., $u$ is a distributional solution of $-\Delta u=\mu$. This follows by choosing $\varphi\in \C_c^\infty(\HS)\subset \D$ in \eqref{eq:defsol}.
An analogue remark is valid for a weak solution of \eqref{dish}.
\end{remark}

\begin{remark}\label{rem:uniq}
We  point out that if $u$ is a weak solution of
\begin{equation*}%\label{eqh1}
 \left\{\begin{array}{ll}
-\Delta u = \mu_1& \mathrm{on }\ \HS, \cr \cr
u= \nu_1 & \mathrm{on }\ \DHS, \cr
\end{array}
 \right.
\end{equation*}
 and of
\begin{equation*}
 \left\{\begin{array}{ll}
-\Delta u = \mu_2& \mathrm{on }\ \HS, \cr \cr
u= \nu_2 & \mathrm{on }\ \DHS, \cr
\end{array}
 \right.
\end{equation*}
where $\mu_1,\mu_2 \in {\cal M}_{loc}(\CHS,x_N)$ and $\nu_1,\nu_2\in {\cal M}_{loc}(\DHS),$ 
%are  local Radon measures on $\Rset^N_+$ and  $\DHS$ respectively,
 then $\mu_1=\mu_2$ and $\nu_1=\nu_2$. Indeed, from the definition of weak solution, we have that for any $\varphi\in  \D$
\begin{equation}\label{eq:double}
\begin{aligned}
  \int_{\HS}\varphi(x) d\mu_1(x)+ \int_{\DHS} \varphi_N(x',0) d\nu_1(x')= 
	\int_{\HS}u (-\Delta \varphi)\qquad \qquad\\ =
  \int_{\HS}\varphi(x) d\mu_2(x)+ \int_{\DHS} \varphi_N(x',0) d\nu_2(x').
\end{aligned}
\end{equation}
%{\color{white} Therefore, by choosing $\varphi\in \C_c^\infty(\HS)\subset \D$, we obtain,
%$$\int_{\HS}\varphi(x) d\mu_1(x)=\int_{\HS}\varphi(x) d\mu_2(x).$$
%Since this is true for any $\varphi\in \C_c^\infty(\HS),$ SI DEDUCE DAL REMARK}
From Remark \ref{rem:weak>distr}
 we deduce that $\mu_1=\mu_2$,
and by \eqref{eq:double} we get,
$$ \int_{\DHS} \varphi_N(x',0) d\nu_1(x')=  \int_{\DHS} \varphi_N(x',0) d\nu_2(x'),\qquad\forall \varphi\in  \D.$$
Next, let $\psi\in \C_c^\infty(\Rset^{N-1}).$ Let  $\phi\in \C^2_c(\Rset)$  be a 
standard cut-off function, that is
\begin{equation}\label{cutoff}
0\le\phi\le 1,\quad \phi(t)=1\ for\ 0\le t\le1, \quad \phi(t)=0\ for\ t\ge 2,
\end{equation}
 %standard cut-off function that is
%\begin{equation}\label{cutoff}
%0\le\phi\le 1,\quad \phi(t)=1\ for\ 0\le t\le1, \quad \phi(t)=0\ for\ t\ge 2,
%\end{equation}
 and set, 
$\varphi(x',x_N)\decl x_N\psi(x')\phi(x_N)$.
Since $\partial_N\varphi(x',0)=\psi(x')\phi(0)=\psi(x')$, we get,
$$ \int_{\DHS} \psi(x') d\nu_1(x')=  \int_{\DHS} \psi(x') d\nu_2(x').$$
Hence $\nu_1=\nu_2$.
\end{remark}

\subsection{Weak solutions and lim-trace}\label{sec3.1}
In this section, we introduce the definition of lim-trace for functions belonging to  $W^{1,p}_{loc}(\HS)$ and we study its properties  in connection to the weak solutions of  \eqref{eqh} and \eqref{dish}.

In this respect, we observe that if $u$ is a weak solution of 
\eqref{eqh} or \eqref{dish} with $\mu$ and $\nu$ local Radon measures, then for any bounded open set 
$\Omega\subset\overline \Omega\subset \HS$, we have $u\in W^{1,p}(\Omega)$ for any $1\le p<N/(N-1)$ (see i.e. \cite{dau-lio80}).
This means that $u$ has a Sobolev trace $u_{\partial\Omega}$ on $\partial\Omega.$ Therefore, if $u$ is a weak solution of 
\eqref{eqh} or \eqref{dish}, then $u\in W^{1,p}_{loc}(\HS)$ and, for any $\epsilon>0$,
$u$ has a Sobolev trace $u(\cdot,\epsilon)$  on each hyperplane with equation $x_N=\epsilon.$ 
Mimicking the argument in \cite{axl} for harmonic functions, we define the trace of a given function on $\DHS$ as the limit of the Sobolev trace $u(\cdot,\epsilon)$ as $\epsilon\to 0^+$.
More precisely, we have the following.
\begin{definition}\label{def:trace} Let
  $u\in W^{1,p}_{loc}(\HS)$ for some $p>1$ and let $\nu$ be a 
  local Radon measure on $\DHS$. We say that $\nu$ is the 
  \emph{lim-trace of $u$ on $\DHS$} and we write $\nu=Tr(u)$, if for any 
  $\psi\in \C^2_c(\DHS)$, we have
  	$$\lim_{\epsilon\to 0^+} \int_{\DHS} \psi(x') u(x',\epsilon) dx' =\int_{\DHS} \psi(x')d\nu(x'),$$
  where $u(\cdot ,\epsilon)$ is the Sobolev trace of $u$ on the
  hyperplane with equation $x_N=\epsilon$.
\end{definition}

Definition \ref{def:trace} can similarly be extended to a function $u\in BV_{loc}(\HS)$.
If $u\in\C(\CHS)$ then $Tr(u)=u(x',0)$. If 
$u\in  W^{1,p}_{loc}(\CHS)$ %\footnote{We write $u\in W^{1,p}_{loc}(\CHS)$ if  $u\in   W^{1,p}(\Omega)$ for any open bounded set $\Omega\subset\HS$.}  
 then
$Tr(u)$ is the Sobolev trace of $u.$ In a more general setting,
if $u\in BV_{loc}(\CHS),$ %\footnote{We write $u\in BV_{loc}(\CHS)$ if  $u\in   BV(\Omega)$ for any open bounded set $\Omega\subset\HS$.}  
 then
$Tr(u)$ is the trace of $u$ in the BV sense 
(see \cite{giu84}). 

We point out that if $u$ and $v$ have lim-traces, then $u+v$ 
also has a lim-trace and $Tr(u+v)=Tr(u)+Tr(v)$. Similarly if $a\in \Rset$,
then $Tr(au)=a Tr(u)$.

The following results are similar to the ones proved in \cite{marver} in another context.
\begin{theorem}\label{teo:traceweaksol} Let $\mu\in {\cal M}_{loc}(\CHS,x_N)$ and
  $\nu\in {\cal M}_{loc}(\DHS)$.
Let $u$ be a weak solution of
\begin{equation}\label{eq:coer}
 \left\{\begin{array}{ll}
 \Delta u = \mu & \mathrm{on }\ \HS, \cr \\
u = \nu  & \mathrm{on }\ \DHS. \cr
\end{array}
 \right.
\end{equation}
Then $u$ has lim-trace on $\DHS$,  and 
$Tr(u)= \nu$.
\end{theorem}

Another important result for functions that possess lim-trace is for nonnegative superharmonic distributions
considered in Corollary \ref{cor:superhamonictrace} below.
\medskip

We need the following Lemma.
\begin{lemma}\label{lem:intpart} Let $\mu\in  {\cal M}_{loc}(\HS)$ and
  let $u$ be a distributional solution of 
  $$ \Delta u = \mu \ \ in\ \ {\cal D}'(\HS).$$
  Then, for any Lipschitz bounded open set 
  $\Omega\subset\overline \Omega\subset \HS$, it follows  that
  $u\in W^{1,p}(\Omega)$ for $1\le p <N/(N-1)$. Moreover, for any $\varphi\in \D(\Omega)$ 
  there holds,
  \begin{equation}\label{eq:intbypart}
    \int_{\partial\Omega} u_{\partial\Omega}(\nabla\varphi\cdot n)=\int_\Omega u\Delta \varphi-\int_\Omega\varphi d\mu.
  \end{equation}
\end{lemma}
\begin{proof}
  Since $u\in W^{1,p}(\Omega)$,
  for any $\varphi\in\C^2_c(\Rset^N),$ we have
  \begin{equation}\label{eq:intbyparttec}
  \int_\Omega u\Delta\varphi=-\int_\Omega \left(\nabla\varphi\cdot\nabla u\right)+ \int_{\partial\Omega} u_{\partial\Omega}(\nabla\varphi\cdot n),
\end{equation}   
(see i.e. \cite{leoni}).
 Next, let $(m_\eta)$ be an approximation of the identity with $supp(m_\eta)\subset B_\eta(0)$,
 and let us choose $\eta$ sufficiently small such that
 a compact $\eta$-neighborhood of $\Omega$ is still contained in $\HS$,
 $\Omega+B_\eta\subset\overline{\Omega+B_\eta}\subset\HS$.
 Since $u_\eta\decl m_\eta \ast u$ is a classical solution of 
 $\Delta u_\eta=\mu_\eta$, where $\mu_\eta\decl  m_\eta\ast\mu$,
 integrating by parts we obtain
  $$\int_\Omega \varphi\mu_\eta=\int_\Omega \varphi\Delta u_\eta=
 \int_{\partial\Omega}\varphi (\nabla u_\eta\cdot n)-
 \int_\Omega \left(\nabla\varphi\cdot\nabla u_\eta\right).$$
  Now, if $\varphi\in  \D(\Omega)$ the first integral in the
  right hand side of the above identity vanishes. 
  Letting $\eta\to 0$ in the above identity, it follows that
	$$ \int_\Omega \varphi d\mu=
	-\int_\Omega \left( \nabla\varphi\cdot\nabla u\right),$$
 which, combined with \eqref{eq:intbyparttec} yields the claim.
\end{proof} 

\begin{proof}[Proof of Theorem \ref{teo:traceweaksol}]  
  Fix  $\psi\in \C^2_c(\DHS)$. For $\epsilon>0,$ define  $\varphi_\epsilon(x',x_N)\decl (x_N-\epsilon)\psi(x')\phi(x_N)$
  where $\phi\in\C^2_c(\Rset)$ is a standard   cut-off function as in \eqref{cutoff}.
  Let $R>0$ be sufficiently large so that $supp(\psi)\subset B'_R\decl \{|x'|<R\}$.
For $0<\epsilon<1$, define the cylinder $\Omega_\epsilon$ be the cylinder $\Omega_\epsilon\decl B'_R\times (\epsilon,2)$.
Since $\varphi_\epsilon\in \D(\Omega_\epsilon)$, applying the integration by parts formula
\eqref{eq:intbypart} gives:
$$ \int_{\Omega_\epsilon} \varphi_\epsilon d\mu - \int_{\Omega_\epsilon} u\Delta \varphi_\epsilon=-\int_{\partial\Omega_\epsilon} u_{\partial\Omega_\epsilon}(\nabla\varphi_\epsilon\cdot n)=\int_{\DHS} u(x',\epsilon)\psi(x')dx',
$$
where we have used the fact that $\nabla \varphi_\epsilon =0$ on 
$\partial\Omega_\epsilon\cap \{x_N> \epsilon\} $ and 
$(\nabla\varphi_\epsilon\cdot n)=-\partial_N\varphi_\epsilon=-\psi$ on $\partial\Omega_\epsilon\cap \{x_N= \epsilon\} $.
Next, observe that as $\epsilon \to 0^+$, $\varphi_\epsilon\to \varphi_0$ uniformly, where $\varphi_0(x',x_N)= x_N\psi(x')\phi(x_N).$
Since $\mu \in {\cal M}_{loc}(\CHS,x_N),$ we have:
$$\int_{\Omega_\epsilon} \varphi_\epsilon d\mu \to \int_{\Omega_0} \varphi_0 d\mu\qquad as\ \ \epsilon\to 0^+.$$
On the other hand, since again $\Delta\varphi_\epsilon \to \Delta\varphi_0$ uniformly, 
and  $u\in L^1(\Omega_0)$, we obtain
$$\int_{\Omega_\epsilon} u\Delta \varphi_\epsilon \to  \int_{\Omega_0} u\Delta \varphi_0 \qquad as\ \ \epsilon\to 0^+.$$
 Noting that $\varphi_0(x',x_N)\decl x_N\psi(x')\phi(x_N)$
 we can conclude that there exists the limit in the definition of lim-trace, that is,
\begin{equation}\label{eq:trace}\begin{aligned}
T(\psi)\decl\lim_{\epsilon\to 0^+} \int_{\DHS} \!u(x',\epsilon)\psi(x')dx'= 
\lim_{\epsilon\to 0^+}\left(  \int_{\Omega_\epsilon}\! \varphi_\epsilon d\mu - \int_{\Omega_\epsilon}\! u\Delta \varphi_\epsilon\right)\\
= \int_{\HS} \!\varphi_0 d\mu - \int_{\HS}\! u\Delta \varphi_0.\qquad\qquad\qquad
\end{aligned}
\end{equation}
To complete the proof, it remains to show that $T$ does not depend on the choice of $\phi$ (which implies that $T$ is a distribution on $\Rset^{N-1}$)
and can be represented by a local Radon measure.
To this end,  since $u$ is a weak solution of equation \eqref{eq:coer}, we utilize
$\varphi_0(x',x_N)\decl x_N\psi(x')\phi(x_N)$ in the definition of weak solution. We deduce that
 $$T(\psi)=\int_{\HS} \varphi_0 d\mu - \int_{\HS} u\Delta \varphi_0=
\int_{\DHS}\partial_N\varphi_0d\nu =\int_{\DHS}\psi d\nu,$$
which concludes the proof.
\end{proof}

The following can be seen as a sort of converse of Theorem \ref{teo:traceweaksol}.
\begin{theorem}\label{teo:traceweakinv} Let $\mu\in {\cal M}_{loc}(\CHS,x_N)$ and
  let $u\in %W^{1,p}_{loc}(\HS)\cap 
L^1_{loc}(\CHS)$ be such that,
$$ \Delta u = \mu \ \ in\ \ {\cal D}'(\HS).$$
If $u$ has lim-trace $\nu=Tr(u)$ with 
$\nu\in {\cal M}_{loc}(\DHS)$, then $u$ is a weak solution of
\eqref{eq:coer}.
\end{theorem}
\begin{remark} A further connection  between superharmonic functions on $\HS$ with a prescribed behavior on the boundary $\DHS$ and weak solutions can be found in Corollary \ref{cor:linkDW} below.
\end{remark}
\begin{proof}[Proof of Theorem \ref{teo:traceweakinv}] We need to demonstrate that $u$ is a weak solution of
\eqref{eq:coer} with $\nu=Tr(u)$. In other words, we must prove that for any $\varphi\in \D$ the following holds:
$$  \int_{\HS}\varphi(x) d\mu(x)
	- \int_{\HS}u \Delta \varphi=\int_{\DHS} \varphi_N(x',0) d\nu(x').$$
Let $\epsilon$ be such that  $0<\epsilon<1$ and define $\tilde{\varphi_\epsilon}\decl\varphi(x',x_N-\epsilon).$
Let $\Omega_\epsilon\decl B'_R\times (\epsilon,R+1)$ with $R$ large enough such that, 
$supp(\varphi)\subset B'_{R}\times (-R,R)$.
Observing that $\tilde\varphi_\epsilon=0$ on 
$\partial\Omega_\epsilon$,
an application of \eqref{eq:intbypart} yields,
$$\begin{aligned}
 \int_{\Omega_\epsilon} \tilde\varphi_\epsilon d\mu - \int_{\Omega_\epsilon} u\Delta \tilde\varphi_\epsilon=&-\int_{\partial\Omega_\epsilon} u_{\partial\Omega_\epsilon}(\nabla\tilde\varphi_\epsilon\cdot n)=\int_{x'\in\Rset^{N-1},x_N=\epsilon} u(x',\epsilon)\partial_N\tilde \varphi_\epsilon dx'\\
=&\int_{\DHS} u(x',\epsilon)\partial_N \varphi(x',0) dx'.
\end{aligned}$$
Taking the limit as $\epsilon\to 0$ we obtain the desired result.
\end{proof}

We conclude this subsection with some generalizations of Theorem \ref{teo:traceweaksol}
and \ref{teo:traceweakinv}.

\begin{lemma}\label{lem:trace}  Let $\mu\in {\cal M}_{loc}(\CHS,x_N)$ and
  $\nu_1\in {\cal M}_{loc}(\DHS)$.
Suppose $u$ is a weak solution of
\begin{equation}\label{dis1:coer}
 \left\{\begin{array}{ll}
 \Delta u = \mu & \mathrm{on }\ \HS, \cr \\
u \le \nu_1  & \mathrm{on }\ \DHS. \cr
\end{array}
 \right.
\end{equation}
Then, $u$ has lim-trace on $\DHS$. Furthermore, $Tr(u)\le \nu_1,$ and $u$ is a weak solution of
\eqref{eq:coer} with $\nu=Tr(u)$.
\end{lemma}
\begin{proof} Following a similar approach to the proof of  Theorem \ref{teo:traceweaksol}, and using the same notations, consider a function $\psi\in\C^2_c(\Rset^{N-1}).$ The functional
$T(\psi),$ defined in \eqref{eq:trace} is a distribution on $\Rset^{N-1}$.
Given that $u$ is a weak solution of \eqref{dis1:coer}, using $\varphi_0$ in the definition of the solution we obtain:
$$\int_{\DHS}\psi d \nu_1\ge\int_{\HS} \varphi_0 d\mu - \int_{\HS} u\Delta \varphi_0=T(\psi).$$
This implies that the expression $$\psi\to \int_{\DHS}\psi d \nu_1 - T(\psi)$$ is a nonnegative distribution on 
$\Rset^{N-1}.$ Consequently, there exists a positive Radon measure $\nu_2$ on $\Rset^{N-1}$,
that represents this distribution. In other words,
$\nu_1-T=\nu_2,$ indicating that $T$ is a local Radon measure. This establishes that $u$ has lim-trace, and $Tr(u)\le \nu_1$.

Moreover, since
$$ \Delta u = \mu \ \ in\ \ {\cal D}'(\HS),$$
 applying Theorem \ref{teo:traceweakinv},
leads to the conclusion that  $u$ is a weak solution of
\eqref{eq:coer} with $\nu=Tr(u)$. 
\end{proof}

\begin{theorem}  Let $\mu_1\in {\cal M}_{loc}(\CHS,x_N)$ and $\nu_1\in {\cal M}_{loc}(\DHS)$.
Suppose $u$ is a weak solution of
\begin{equation}\label{dis:coer1}
 \left\{\begin{array}{ll}
 \Delta u \ge \mu_1 & \mathrm{on }\ \HS, \cr \\
u \le \nu_1 & \mathrm{on }\ \DHS. \cr
\end{array}
 \right.
\end{equation}

1. Assume that for any $R>0,$ there exist $M>0$ $\delta>0$ such that
$u\le M$ on the cylinder $\{|x'|<R\}\times (0,\delta)$. Then $u$ has lim-trace on  $\DHS$.

2. Assume that for any $R>0,$ there exist $\tilde M>0$ $\delta>0$ such that
$\int_{|x'|<R}|u(x',\epsilon)| dx' <\tilde M$ for any $0<\epsilon<\delta$. 
Then $u$ has lim-trace on  $\DHS$.

3. Assume that $u$ has lim-trace on  $\DHS$. Then $Tr(u)\le \nu_1,$ and
there exist a positive Radon measure $\lambda\in {\cal M}_{loc}(\CHS,x_N)$ such that $u$ is a solution of 
\eqref{eq:coer} with $\mu=\mu_1+\lambda$ and $\nu=Tr(u)$.
\end{theorem}
\begin{proof}
  Since $u$ satisfies 
  $$\Delta u\ge \mu_1 \quad in \ {\cal D}'(\HS),$$
  there exists a positive radon measure $\lambda\in {\cal M}^+(\HS)$ such that
  $u$ solves $$\Delta u = \mu_1+\lambda \quad in \ {\cal D}'(\HS).$$

  Arguing as in the proof of Theorem \ref{teo:traceweaksol}, and adopting the same notations, 
  for a nonnegative function $\psi\in\C^2_c(\Rset^{N-1}),$ we have
  $$I^\epsilon_1:= \int_{\DHS} u(x',\epsilon)\psi(x')dx'=
   \int_{\Omega_\epsilon} \varphi_\epsilon d\mu_1 - \int_{\Omega_\epsilon} u\Delta \varphi_\epsilon +   \int_{\Omega_\epsilon} \varphi_\epsilon d\lambda
   =:I^\epsilon_2-I^\epsilon_3+I^\epsilon_4.$$
Next,  as in \eqref{eq:trace}
  $I^\epsilon_2-I^\epsilon_3$ converges as $\epsilon \to 0$. 
  By the Beppo Levi monotone convergence theorem, since 
  $\varphi_\epsilon\chi_{\Omega_\epsilon}\nearrow\varphi_0$, it follows that
  $I^\epsilon_4\to  \int_{\HS} \varphi_0 d\lambda$, so $I^\epsilon_1$
  has a finite or infinite limit. 
  We claim that,
\begin{equation}\label{claim0}
 \int_{\HS} \varphi_0 d\lambda<\infty.
\end{equation}  
This  claim implies that  $\lambda\in {\cal M}_{loc}(\CHS,x_N)$.
Note that \eqref{claim0} is equivalent to  $I^\epsilon_1$ having a finite limit
as $\epsilon\to 0$.
  
Now, suppose that 1 holds. Let $R>0$ be large enough  such that $supp(\psi)\subset B'_R$.
Then, by our assumption there exist $M,\delta>0$ such that  for $\epsilon<\delta,$ we have
  $I^\epsilon_1\le M \int_{\DHS} \psi(x')dx'.$ Therefore, 
  $I^\epsilon_4$  is uniformly bounded for $\epsilon<\delta$.
  This implies that $I^\epsilon_1$ has a finite limit as $\epsilon\to 0$. That is,
  $T(\psi)\decl\lim_{\epsilon\to 0^+} \int_{\DHS} u(x',\epsilon)\psi(x')dx'$ exists and is finite.
Using a similar argument as in the proofs of  Theorem \ref{teo:traceweaksol} and Lemma \ref{lem:trace},
  we conclude that $T$ can be represented by a local Radon measure. This completes the proof. 
  
  Next, suppose that 2 holds. The claim will follow immediately by arguing as in the case 1, using the estimate
 $| I^\epsilon_1|\le \abs{\abs{\psi}}_\infty\int_{|x'|<R}|u(x',\epsilon)| dx' < 
  \abs{\abs{\psi}}_\infty \tilde M$.
  
Finally, suppose that 3 holds and let $\nu=Tr(u)$. We know that 
  $I^\epsilon_1\to \int_{\DHS} \psi(x')d\nu(x')$ as $\epsilon\to 0$. This implies that,
  the $I^\epsilon_4$ converges to a finite limit, so \eqref{claim0} holds. Hence,
  $\lambda\in {\cal M}_{loc}(\CHS,x_N)$.
  
 The fact that  $u$ is a solution of
\eqref{eq:coer} with $\mu=\mu_1+\lambda$ and $\nu=Tr(u),$ 
  is a consequence of Theorem \ref{teo:traceweakinv}.
\end{proof}
By  slightly modifying the above proofs, we can establish the following generalization of Theorem \ref{teo:traceweakinv}.
\begin{theorem}\label{teo:traceweakinvdis} Let $\mu\in {\cal M}_{loc}(\CHS,x_N),$ and let $u\in  L^1_{loc}(\CHS)$ satisfy
\begin{equation}\label{dis:distr}
 \Delta u \ge \mu_1 \ \ in\ \ {\cal D}'(\HS).
\end{equation}
If $u$ has lim-trace $\nu=Tr(u)$ with 
$\nu\in {\cal M}_{loc}(\DHS)$, then there exist a positive Radon measure $\lambda\in {\cal M}_{loc}(\CHS,x_N)$ such that $u$ is a solution of 
\eqref{eq:coer} with $\mu=\mu_1+\lambda$ and $\nu=Tr(u)$.

In particular, $u$ is a weak solution of
\begin{equation*}
 \left\{\begin{array}{ll}
 \Delta u \ge \mu_1 & \mathrm{on }\ \HS, \cr \\
u= \nu & \mathrm{on }\ \DHS. \cr
\end{array}
 \right.
\end{equation*}
Furthermore, if equality holds in \eqref{dis:distr}, then $\lambda\equiv 0$.
\end{theorem}

For the lim-trace of a weak solution of \eqref{dish} with positive Radon measure, refer to Theorem~\ref{teo:distr=weak} and Corollary \ref{cor:linkDW} provided below.

\begin{remark}[on Lim-Trace of Positive Parts of Functions]\label{rem:counterexample}   It is interesting to investigate whether 
 a function $u$ that has a lim-trace also ensures that its  
 positive part $u^+$ possesses the same property.
Generally, the answer is negative. Consider the following example: Let
  $m_1\in \C^\infty_c(\Rset)$ be an even, nonnegative standard cut off function supported in $(-1,1)$ with $\int m_1=1.$ For $\epsilon>0,$ define the mollifier family as
  $m_\epsilon(x)\decl \frac{1}{\epsilon} m_1(\frac{x}{\epsilon}).$ 
  Next, for $x\in\Rset$ and $\epsilon>0,$ define $u(x,\epsilon)$ as  (see Figure \ref{fig2})
  $$u(x,\epsilon)\decl \frac{1}{\sqrt{\epsilon}}m_\epsilon(x-\epsilon)-
\frac{1}{\sqrt{\epsilon}}m_\epsilon(x+\epsilon).$$

Clearly, $u\in\C^\infty(\Rset^2_+),$ 
 and by the Lagrange theorem, for any  $\psi\in\C^2_c(\Rset),$ we have the estimate:
  $$ \abs{\int_\Rset u(x,\epsilon)\psi(x)}\le 2 \sqrt{\epsilon}||{\psi'}||_\infty.$$
Thus, $u$ has lim-trace $\nu=Tr(u)=0,$ the trivial measure on $\partial\Rset^2_+=\Rset$.

 However, the positive part $u^+(x,\epsilon)= \frac{1}{\sqrt{\epsilon}}m_\epsilon(x-\epsilon)$ does not admit a lim-trace. For a  nonnegative
  $\psi\in\C^2_c(\Rset)$ with $\psi(0)>0$ we have:
 $$\int_\Rset u^+(x,\epsilon)\psi(x) dx\ge \frac{1}{\sqrt{\epsilon}} \,\frac{\psi(0)}{4}\ \ 
 \ \ \ \mathit{for\ \ sufficiently\ \  small\ \ \epsilon},$$
 and hence,
 
 $$\lim_{\epsilon \to 0}\int_\Rset u^+(x,\epsilon)\psi(x) dx\ge \lim_{\epsilon \to 0}\frac{1}{\sqrt{\epsilon}} \,\frac{\psi(0)}{4}=\infty.$$
 \begin{figure}\centering
\includegraphics[scale=.6]{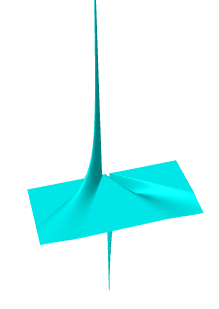}
\caption{The plot of function $u$ of Remark \ref{rem:counterexample}.}\label{fig2}
\end{figure}

On the other hand, in the general case, if both $u$ and  $u^+$ have lim-traces,
 say $\nu=Tr(u)$ and $\lambda=Tr(u^+)$, it can be shown that 
 $\lambda\ge \nu^+$.
 However, it is not necessarily the case that $\lambda=\nu^+.$
 To illustrate this, consider:

  $$v(x,\epsilon)\decl m_\epsilon(x-\epsilon)- m_\epsilon(x+\epsilon).$$
  As before, $v$ has lim-trace and $Tr(v)=0$, while
  $v^+$ admits a lim-trace  $Tr(v^+)=\delta_0,$ the Dirac measure at the origin.
  
The example related to the function $u$ highlights the challenges in relaxing the definition of lim-trace. While lim-trace can be viewed as 
a weak* convergence of the traces $u(x',\epsilon),$ this convergence is achieved through smooth test functions, and the example indicates that it cannot be relaxed to consider only continuous functions. 
Nonetheless, our definition of lim-trace is suitable for our purposes, as
Theorems \ref{teo:traceweaksol} and \ref{teo:traceweakinv} establish a strong connection between 
the weak solutions of  \eqref{eq:coer} and their  lim-traces.
\end{remark}

\section{Main results}\label{sec:main}
To maintain clarity and streamline the organization of this paper, we will postpone the proof of our first main result to Sections \ref{proof1A} and \ref{proof1B}. The reasons for this will become apparent  in due course.

For the reader convenience we rewrite Theorems A and B of the introduction in the following form.

\begin{theorem}\label{teo:repiff}Let $N\geq 2$, and let $\mu$ and $\nu$ be 
 positive Radon measures on $\Rset^N_+$ and  $\DHS,$ respectively.
\begin{enumerate}
\item\label{mainA}
 Let $\mu\in {\cal M}_{loc}(\CHS,x_N)$ and let $u\in L^1_{loc}(\CHS)$  be a weak solution of (\ref{dish}).
Suppose $x\in\HS$ is a Lebesgue point for $u$. Assume that there exists $h\in\Rset$ such that  \rhh holds at the point $x$ that is,
\begin{equation*}
 \liminf_{R\rightarrow +\infty} \frac{1}{R^{N+2}}\int_{\{y_N>0\} \cap \{R<\abss{x-y}<2R\}} y_N |u(y)-h y_N| \ dy =0.
\end{equation*}
   Then,
   \begin{equation}\label{disrepmeas+}
      u(x)\ge hx_N+\int _{\DHS} K^x(y') d\nu(y') +  \int_{\HS} G^x(y) d\mu(y).
   \end{equation}
   
   Moreover, if $u$ is a weak solution of \eqref{eqh}, then the inequality in \eqref{disrepmeas+} becomes  an equality, that is
     \begin{equation}\label{eqrepmeas+}
      u(x)= hx_N+\int _{\DHS} K^x(y') d\nu(y') +  \int_{\HS} G^x(y) d\mu(y).
   \end{equation}
%\begin{enumerate}[label=\arabic*)]
%\item \label{mainAa}
% Let $u\in L^1_{loc}(\CHS)$  be a weak solution of (\ref{eqh}) 
%   and there exists $h\in\Rset$ such that the condition    \rhh holds.
%   Then
%   for a.e. $x\in \Rset^N_+$, we have
%   \begin{equation}\label{eqrepmeas+}
%      u(x)= hx_N+\int _{\DHS} K^x(y') d\nu(y') +  \int_{\HS} G^x(y) d\mu(y).
%   \end{equation}
% \item\label{mainAb}  Let $u\in L^1_{loc}(\CHS)$  be a weak solution of (\ref{dish}) 
%   and there exists $h\in\Rset$ such that the condition   \rhh\ holds.
%   Then  for a.e. $x\in \Rset^N_+$, we have
%   \begin{equation}\label{disrepmeas+}
%      u(x)\ge hx_N+\int _{\DHS} K^x(y') d\nu(y') +  \int_{\HS} G^x(y) d\mu(y).
%   \end{equation}
%\end{enumerate}
\item\label{teo:inverse}
  Let $u$ be defined by \eqref{eqrepmeas+} and assume that it is finite  for a.e. $x\in\HS$.

  Then $u\in L^1_{loc}(\CHS)$, $\mu\in {\cal M}_{loc}(\CHS,x_N)$ and $u$ is a weak solution of \eqref{eqh}.  Moreover, $u$ satisfies
  the condition \rhh\!\!, where the constant $h$ appears in \eqref{eqrepmeas+}.
  
  Furthermore, 
  \begin{enumerate}
     \item $$h=\inf_{\HS} \frac{u(x)}{x_N}.$$
  \item Let $\Omega\subset\HS$ be a nonempty bounded open set. If $h=\inf_{\Omega} \frac{u(x)}{x_N}$,
  then $u(x)=h x_N$, $\mu\equiv 0$ and $\nu\equiv 0$. 
  \item If $u-hx_N$ is not identically zero, then there exists a constant $c_0>0$ such that
  \begin{equation}\label{est:ubelow}
  u(x)\ge hx_N+ c_0 \frac{x_N}{1+|x|^N},\quad for \ \ all\ \ x\in\HS.
  \end{equation}
  \item For almost every $x\in\HS,$
  \begin{equation*}
 \lim_{R\rightarrow +\infty} \frac{1}{R^{N+2}}\int_{\{y_N>0\} \cap B^*_R(x)} y_N |u(y)-h y_N| \ dy =0.
\end{equation*}
  \end{enumerate}
\end{enumerate}
\end{theorem}

Note that the main statement in \ref{teo:inverse}  is the converse of  the statement in \ref{mainA}.
\medskip

As a consequence of Theorems 
 \ref{th:RepHSClassic}, \ref{teo:repiff}.\ref{teo:inverse} and \ref{teo:traceweaksol}, the following corollary holds:  
\begin{corollary}\label{cor:superhamonictrace} A nonnegative superharmonic distribution on $\HS$ possess lim-trace.
\end{corollary}

\begin{remark} We observe that  Theorem \ref{teo:repiff} and Remark \ref{rem:uniq}, imply  the uniqueness of the measures $\mu$ and $\nu$ in the representation Theorem \ref{th:RepHSClassic}.
\end{remark}

\begin{corollary}
 Let $\mu\in L^1_{loc}(\Rset^N_+,x_N)$ and $\nu\in L^1_{loc}(\Rset^{N-1})$ be nonnegative functions, and  let $u$ be a
 weak solution of (\ref{dish})
 satisfying \rhz\!\!. Then for almost every $x\in \Rset^N_+,$ we have
 \begin{equation}
	u(x)\geq \int _{\Rset^{N-1}} K^x(y') \nu(y')dy' +  \int_{\Rset^N_+} G^x(y) \mu(y)dy.
 \end{equation}
\end{corollary}

\begin{corollary}
 Under the assumptions of Theorem \ref{teo:repiff}.\ref{mainA} with $h\ge 0$, if $u$ is a weak solution of (\ref{eqh}), then
 $u(x)>hx_N$ or $u(x)\equiv hx_N$ a.e. on $\HS$.
\end{corollary}

\begin{remark}
In general, a superharmonic function on $\HS$ does not belong to
$L^1_{loc}(\CHS)$. For example, the function $u(x_1,x_2)=-|x|^{-2}$ does not belong to 
$L^1_{loc}(\CHS)$ with $N=2,$ and $-\Delta u= 4 |x|^{-4}\ge 0$. 
See also the example defined in \eqref{eq:unotl1c}.
\end{remark}

\begin{remark} i) In Theorem \ref{teo:repiff}.\ref{mainA},
the fact that,
 the integrals in (\ref{eqrepmeas+}) are well defined is an outcome of our approach.
In other words, the finiteness of the integrals  in (\ref{eqrepmeas+})
 is a necessary condition for the existence of a solution satisfying \rhh\!\!.
 
ii) If we know a priori that the integrals  in (\ref{eqrepmeas+}) are finite, then \eqref{eqrepmeas+} describes all the solutions satisfying \rhh\!\!.

 iii)  Notice  that if we require a priori that the integrals  in (\ref{eqrepmeas+}) are finite without any assumption on the sign of the measures $\mu$ and $\nu$,
 then \eqref{eqrepmeas+} defines a solution of \eqref{eqh}.
 More precisely, if 
 $\int _{\DHS} K^x(y') d|\nu|(y')<\infty$ and $\int_{\HS} G^x(y) d|\mu|(y)<\infty$, then \eqref{eqrepmeas+} defines a weak solution of \eqref{eqh}. For a proof of this, the interested reader can follow the same steps in the proof of Theorem \ref{teo:repiff}.\ref{teo:inverse},
 see Section \ref{proof1B}.
 
A more general result  is given by the following.
\end{remark}

\begin{theorem}\label{repmeasnosign}
 Let $N\geq 2$, and let $\mu\in {\cal M}_{loc}(\CHS,x_N)$ and
  $\nu\in {\cal M}_{loc}(\DHS).$
  %. $\mu$ and $\nu$ be local Radon measures on $\Rset^N_+$ and  $\DHS$ respectively, 

Writing $\mu=\mu^+-\mu^-$ and $\nu=\nu^+ - \nu^-$, for a.e. $x\in\HS$ we assume that
\begin{equation}\label{hypneg}
\int _{\Rset^{N-1}} K^x(y') d\nu^-(y') <\infty,
\quad \int_{\Rset^N_+} G^x(y) d\mu^-(y)<\infty.
  \end{equation} 

 Let $u\in L^1_{loc}(\CHS)$ be a weak solution of (\ref{eqh}).
 If $u$ satisfies \rhh\!\!,  then
 \eqref{eqrepmeas+} holds for a.e. $x\in \Rset^N_+$.
\end{theorem}
The proof of Theorem \ref{repmeasnosign} is a small modification of the 
proof of Theorem~\ref{teo:repiff}.\ref{mainA}, which can be found in 
Section \ref{sec:proof2}.

\begin{remark}\label{rem:ringlim} If $u$ is given by  \eqref{eqrepmeas+}, or equivalently $u$ is a nonegative superharmonic function, then $u$ satisfies \rhh\!\!.
 %as claimed in Theorem \ref{teo:repiff} and 
 Moreover, the $\liminf$ in the \rhh condition is actually a limit, meaning $u$ satisfies
\begin{equation}\label{ringlim}
\begin{aligned}
 there\ exists\ h\in\Rset\ such\ that\qquad\qquad\qquad\qquad \\ 
 \lim_{R\rightarrow +\infty} \frac{1}{R^{N+2}}\int_{\{y_N>0\} \cap \{R<\abss{x-y}<2R\}} y_N |u(y)-h y_N| \ dy =0,\  for\ a.e.\ x\in\Rset^N_+.
\end{aligned}
\end{equation}

Generally, the functional class satisfying \rhh is not a linear space. However, the set of functions that fulfill \eqref{ringlim} forms a linear space.
\end{remark}

\begin{theorem}
Let %$N\geq 2$, and let %$\mu\in {\cal M}_{loc}(\CHS)$ 
$\mu$ be a positive Radon measure such that
$\mu\not\in {\cal M}_{loc}(\CHS,x_N)$.
Then the problem 
$$ -\Delta u=\mu,\quad in\ \ {\cal D}'(\HS),\qquad u\ge 0, $$
does not admit a solution.
\end{theorem}
The proof is an immediate consequence of Theorem \ref{teo:repiff}.\ref{teo:inverse}
(indeed if the problem admits a solution, then necessarily $\mu\in {\cal M}_{loc}(\CHS,x_N)$).

\medskip

 The following result establishes a connection between the notion of weak solution and distributional solution of the problem:
\begin{equation}\label{eq:superharmonic}
  -\Delta u=\mu,\quad in\ \ {\cal D}'(\HS).
 \end{equation}
 An additional result in this context is presented in  Corollary~\ref{cor:linkDW} below.

 Let's first recall the Riesz decomposition theorem (see \cite[Theorem 4.4.1]{AG});
 if $u$ is a superharmonic function satisfying
\eqref{eq:superharmonic} and there exists a superharmonic function $g$ such that
$u\ge -g$, then there exists an harmonic function $H$ on $\HS$, such that
$$u(x)=H(x)+ \int_{\HS} G^x(y) d\mu(y).$$
Moreover, $H$ is the greatest harmonic minorant of $u$.

We shall prove that this harmonic function $H$ can be obtained as a weighted limit of $u$ on certain rings,
see \eqref{eq:ringh} below.
This construction is further discussed in Section \ref{sec:ringreen}, see Theorem \ref{teo:repold}.

The key point in order to establish a connection between a solution $u$ of \eqref{eq:superharmonic} and a weak solution of \eqref{eqh} is whether or not $u$ possess lim-trace.
If $u\ge-g$ with $g$ a superharmonic function, then $u+g$ has lim-trace as stated in Corollary \ref{cor:superhamonictrace}. However, this does not necessarily imply that $u$ itself has a lim-trace. A sufficient condition that ensures $u$ has lim-trace is if $g\ge 0$.

\begin{theorem}\label{teo:distr=weak}
 Let $\mu$ be a positive Radon measure on $\HS,$ and let $u\in L^1_{loc}(\HS)$ be a distributional solution of  \eqref{eq:superharmonic}.

 Assume there exists a nonnegative superharmonic function $g$ on $\HS$
  such that $u\ge - g$ a.e. on $\HS$. 
  Then $u$ admits lim-trace $Tr(u),$ which is a local Radon measure (not necessarily positive) and $\mu\in {\cal M}_{loc}(\CHS,x_N)$. Moreover,  setting $\nu\decl Tr(u) \in {\cal M}_{loc}(\DHS),$ we have:

  (a)   $u$ is a weak solution of  \eqref{eqh}, 
  
  (b) There exists $h_u\in\Rset$ such that
  \eqref{eqrepmeas+} holds a.e. on $\HS$ with $h=h_u$, i.e.
 \begin{equation}\label{repaux}
      u(x)= h_ux_N+\int _{\DHS} K^x(y') d\nu(y') +  \int_{\HS} G^x(y) d\mu(y).
   \end{equation}  
  
  (c)  $u+g$ satisfies \rhh for a suitable $h\ge 0$.
  
  (d)  $u$ satisfies \rhh with $h=h_u$.
  
  (e) If $u$ is bounded from below, then $h_u\ge 0$.
  
   (f)  If $u$ is nonnegative, then the lim-trace $Tr(u)$ is a positive Radon measure.
\end{theorem}
\begin{proof} The first step relies on the same reasoning as in Corollary \ref{cor:superhamonictrace}. Indeed, since $g$ is nonnegative and superharmonic,
  Theorem \ref{th:RepHSClassic} applies (see also Theorem 1.37 in \cite{AG}), implying that there exists
 $h_g\ge 0$, a positive Radon measure $\nu_g$ on $\DHS$ and a positive 
 Radon measure $\mu_g$ on $\HS$, such that
$$g(x)=h_g x_N + \int _{\DHS} K^x(y') d\nu_g(y')+  \int_{\HS} G^x(y) d\mu_g(y),\qquad for\ a.e.\ x\in\HS.$$
Thus, by Theorem \ref{teo:repiff}.\ref{teo:inverse}, it follows that $\mu_g\in {\cal M}_{loc}(\CHS,x_N)$ and 
$g$ is a weak solution of
\begin{equation}\label{pr:har}
 \left\{\begin{array}{ll}
-\Delta g = \mu_g& \mathrm{on }\ \HS, \cr 
g= \nu_g & \mathrm{on }\ \DHS. \cr
\end{array}
 \right.
\end{equation}
Therefore, by Theorem \ref{teo:traceweaksol} 
$g$ has lim-trace $Tr(g)=\nu_g$.

Next, let $v\decl u+g$. 
The function $v$ is nonnegative and satisfies $-\Delta v=\mu+\mu_g\ge0$ in distributional sense on $\HS$.
Using the same argument as above, there exists $h_v\ge 0$ and 
positive Radon measures $\nu_v$ on $\DHS$, $\mu_v\in {\cal M}_{loc}(\CHS,x_N)$ such that 
\begin{equation}\label{repv}
  v(x)=u(x)+g(x)=h_v x_N +\int _{\DHS} K^x(y') d\nu_v(y') +  \int_{\HS} G^x(y) d\mu_v(y).\end{equation}
Additionally,  $v$ is a weak solution of
\begin{equation}\label{pr:aux}
 \left\{\begin{array}{ll}
-\Delta v = \mu+\mu_g& \mathrm{on }\ \HS, \cr 
v= \nu_v & \mathrm{on }\ \DHS. \cr
\end{array}
 \right.
\end{equation}
Hence,  $\mu_v=\mu+\mu_g$ and $v$ has lim trace $Tr(v)=\nu_v$.
Consequently,  $u=v-g$ has lim-trace $Tr(u)=Tr(v)-Tr(g)=\nu_v-\nu_g\in {\cal M}_{loc}(\DHS)$ and $\mu\in {\cal M}_{loc}(\CHS,x_N)$.

By applying linearity, from \eqref{pr:har} and \eqref{pr:aux}, we complete the proof of (a).

The proof of (b) follows by reformulating the integral expressions of $v$ and $g$
mentioned earlier, where  $h_u=h_v-h_g$.
The claim of (c) is a direct result of applying Theorem \ref{teo:repiff}.\ref{teo:inverse} to \eqref{repv}.

 Proof of (d). Utilizing Theorem \ref{teo:repiff}.\ref{teo:inverse} and Remark \ref{rem:ringlim} applied to $g$ and $v$, it can be deduced that both $g$ and $v$ satisfy \eqref{ringlim}. Consequently, $u=v-g$ also satisfies \eqref{ringlim} and therefore \rhh\!\!.
 
 Proof of (e). If $u$ is bounded from below, i.e., if $g$ is a positive constant  $g=a$, then it follows that $\nu_g=a.$ In other words, $\nu_g$ is a positive multiple of the Lebesgue measure on $\DHS,$ and $h_g=0$. Hence, $\nu_u=\nu_v -a,$ and 
$h_u=h_v-0\ge 0$.

 Proof of (f). This conclusion follows directly from the previous points.\end{proof}

The following result provides a characterization of a comparison principle.
\begin{theorem}\label{teo:potentialink} 
Let $u,v\in L^1_{loc}(\HS)$ such that
$$-\Delta u\ge -\Delta v \ \ in\ \  {\cal D}'(\HS),$$
and 
\begin{equation}\label{liminfu-v}
\liminf_{x\to (y',0))} u(x)\ge \limsup_{x\to (y',0)} v(x),\quad for\ any\ \ y'\in\DHS.\end{equation}
Then $x_N|u-v|\in L^1_{loc}(\CHS)$.

Furthermore,
\begin{equation}\label{umv}
 u\ge v\quad a.e.\ on\ \  \HS,
\end{equation} if and only if 
\begin{equation}\label{u-vRing}
 u-v\  \ satisfies\ \ ({\cal R}^+)\  \ with\  h\ge 0.
\end{equation} 
Moreover, if \eqref{umv}, or equivalently \eqref{u-vRing}, holds, then 
  $u-v\in L^1_{loc}(\CHS)$.
\end{theorem}

In particular, this comparison principle applies to functions 
that satisfy \rhz\!\!.
For another comparison principle, see Section 4 in \cite{damiII}.

To avoid unnecessary complexity, we defer the proof
of Theorem~\ref{teo:potentialink} to Section \ref{secproofpotlink}. The reason is that it relies on notations and results introduced the proof 
of Theorem~\ref{teo:repiff}.

\begin{remark}\label{reml1loc}
In the above theorem, the fact that $u-v\in L^1_{loc}(\CHS)$ is a significant consequence of  \eqref{umv} or \eqref{u-vRing}. 
%Indeed, if we know that $u-v$ is locally integrable up to the boundary,
%the proof of the implication $\eqref{u-vRing}\Rightarrow \eqref{umv}$ can be simplified.

Furthermore, consider $N\ge 2$ and $v\equiv 0$. Define
\begin{equation}\label{eq:unotl1c}
u(x)\decl \frac{1}{|x|^{N}}\left(1-N\frac{x_N^2}{|x|^2}\right).\end{equation}
The function $u$ is harmonic on $\HS$ and satisfies \rhz\!\!.  However, it is not nonnegative, which is  because $u\not\in  L^1_{loc}(\CHS)$. 
In this scenario, the conditions of Theorem \ref{teo:potentialink} are not satisfied.
 Indeed, if $y'\neq 0$, we have  $\lim_{x\to (y',0)} u(x)=|y'|^{-N}>0$, while
if $y'=0$ we have $\liminf_{x\to (0,0)} u(x)=-\infty$.
This implies that  condition \eqref{liminfu-v} must holds for every $y'\in\DHS$. 
The latter cannot be relaxed to hold for almost every  $y'\in\DHS$.
\end{remark}

As a consequence of Theorems \ref{teo:potentialink} and \ref{teo:repiff}, 
we obtain the following result, which establishes a connection between the notion of weak solution and distributional solution.

\begin{corollary}\label{cor:linkDW} Let $u\in L^1_{loc}(\HS)$ be a superharmonic function such that
\begin{equation}\label{eq:boundary}
\liminf_{x\to (y',0)} u(x)\ge 0,\quad for\ any\ \ y'\in\DHS.\end{equation}
Then $x_N|u|\in L^1_{loc}(\CHS)$.

Furthermore, $u\ge 0$ if and only if $u$ satisfies \rhh with $ h\ge 0$.

Moreover, if  $u$ satisfies \rhh for some $h\in\Rset$, then  there exist suitable positive Radon measures $\mu$ and $\nu$ such that
 \eqref{eqrepmeas+} holds and $u$ is a weak solution of \eqref{eqh}.
\end{corollary}

Related results on the maximum principle in bounded domains
for superharmonic function with boundary condition like in \eqref{eq:boundary} can be found in \cite{BNV}.
\medskip

The following positivity result seems noteworthy on its own.
\begin{theorem}\label{teo:bb=0}  Let $\mu$ and $\nu$ be positive
 Radon measures on $\Rset^N_+,$ and  $\DHS,$ respectively. Let $u$ be a weak solution of \eqref{dish} [or \eqref{eqh}].
  If $u$ is bounded from below, then $u$ is nonnegative.  
  Furthermore, $u$ has lim-trace, satisfies \rhh with $h\ge 0$, can be represented by
  \eqref{disrepmeas+}  [or \eqref{eqrepmeas+}], and Theorem \ref{teo:repiff} applies.
\end{theorem}
\begin{proof} 
  Since $u$ is superharmonic on $\HS$, Theorem \ref{teo:distr=weak} is applicable.
  From points (d) and (e), it follows that $u$  satisfies \rhh with $h\ge 0$. Therefore, applying Theorem \ref{teo:repiff} completes the proof.
\end{proof}

\subsection{Proof of Theorem \ref{teo:repiff}.\ref{mainA}}\label{proof1A}

%\begin{lemma}Let $N\geq 2$, let $\mu$ and $\nu$ be 
% positive Radon measures on $\Rset^N_+$ and  $\DHS$ respectively.
% Let $u\in L^1_{loc}(\CHS)$  be a weak solution of (\ref{dish}).
% Let $x\in\HS$ be a Lebesgue point for $u$ and assume that
% and there exists $h\in\Rset$ such that  the condition  \rhh\ holds at the point $x$, that is
%\begin{equation*}
% \liminf_{R\rightarrow +\infty} \frac{1}{R^{N+2}}\int_{\{y_N>0\} \cap \{R<\abss{x-y}<2R\}} y_N |u(y)-h y_N| \ dy =0.
%\end{equation*}
%   Then
%   \begin{equation}\label{replem}
%      u(x)\ge hx_N+\int _{\DHS} K^x(y') d\nu(y') +  \int_{\HS} G^x(y) d\mu(y).
%   \end{equation}
%   
%   Moreover, if $u$ is a weak solution of \eqref{eqh}, then \eqref{replem} holds with the equality sign.
%\end{lemma}
The proof employs a systematic method to demonstrate the claim by utilizing Green's functions and test functions in a weak formulation that follows the following breakdown:
\begin{proof}
%Step 1: Adjusting the Function $u.$ Initially, 
We adjust $u$ by replacing it with $u(y)-hy_N$
to ensure that the function satisfies the condition \rhz at the point $x.$
This is a necessary step to simplify the proof.
%Let $x\in \HS$ be such that \rhz holds.

In what follows, we use the notation of Appendix \ref{sec2} along with the results contained therein.
%Step 2: Choice of Test Function. 
We select the test function
$\tilde \varphi(y)\decl G_\epsilon^x(y) \varphi(y),$ where $G_\epsilon^x$ is the regularized Green's function centered at $x$ with parameter $\epsilon>0$ (see \eqref{Geps}), and $\varphi\in \C_c^2(\CHS)$.
% is a compactly supported smooth function. The Green's function is smooth and satisfies the condition $G_{\epsilon}^x(y',0)=0$, allowing us to apply it within the weak solution definition.
Since $G^x_\epsilon\in C^\infty$ and 
\begin{equation}\label{G0}
G^x_\epsilon(y',0)=0,
\end{equation}
 we can use
$\tilde \varphi$ as test function in \eqref{eq:defsol} obtaining,
\begin{equation}\label{repepsilon}
  \int_{\HS}G^x_\epsilon(y)\varphi(y) d\mu(y)+ 
  \int_{\DHS} K^x_\epsilon(y')\varphi(y',0)  d\nu(y')\le 
	-\int_{\HS}u(y) \Delta( G^x_\epsilon \varphi)(y) dy\ 
\footnote{Here, we have used the definition of $K^x$, \eqref{G0}
  and the fact that,
$$  \frac{\partial}{\partial y_N}(G^x_\epsilon \varphi)(y',0)=
  \frac{\partial G^x_\epsilon}{\partial y_N}(y',0) \varphi (y',0)+ G^x_\epsilon(y',0) \frac{\partial\varphi}{\partial y_N}(y',0)= K^x_\epsilon(y')\varphi(y',0).
$$ 
}\end{equation}
$$\begin{aligned}
%= -\int_{\HS}u(y) \Delta( G^x_\epsilon \varphi)(y) dy\\
=\int_{\HS} u(y)\varphi(y) (-\Delta G^x_\epsilon)(y) dy+
\int_{\HS}u(y) (-\Delta \varphi)(y) G^x_\epsilon(y) dy\\
-2 \int_{\HS}u(y)\nabla G^x_\epsilon\cdot\nabla\varphi=: I_1+I_2+I_3.
\end{aligned}
$$
For the first integral we have,
$$\begin{aligned}
I_1=&\int_{\HS} u(y)\varphi(y) (-\Delta G^x_\epsilon)(y) dy\\
=&
\int_{\HS} u(y)\varphi(y)   (-\Delta \Gamma^x_\epsilon)(y) dy-
\int_{\HS} u(y)\varphi(y)   (-\Delta \Gamma^{\hat x}_\epsilon)(y) dy=:I_{11}+I_{12}.
\end{aligned}
$$
 Now from the estimate,
$$|u(y)\varphi(y)  (-\Delta \Gamma^{\hat x}_\epsilon)(y)|=
|u(y)\varphi(y)  c \frac{\epsilon^2}{\xg^{(N+2)/2}}|\le 
|u(y)\varphi(y)| c \frac{\epsilon^2}{x_N^{N+2}}\in L^1(\HS),$$
an application of the Lebesgue dominated convergence theorem implies that $I_{12}\to 0$ as $\epsilon\to 0$.

In addition, since $x$ is a Lebesgue point of $u,$
$I_{11}\to u(x)\varphi(x)$  as $\epsilon\to 0$.

Finally, in order to estimate  $I_2$ and $I_3$ we choose $\varphi$ as 
$$\varphi(y)\decl\phi(\frac{|x'-y'|}{R})\phi(\frac{|x_N-y_N|}{R}),$$
where $R>0$ and $\phi\in \C^2_c(\Rset)$ is a standard cut-off function as in \eqref{cutoff}.

Note that the support of $\varphi$ is contained in 
$\{y\in\Rset^N: |y-x|_*\le 2R\}$,
and $\varphi\equiv 1$ on $B^*_{R}(x)$,
while
the support of $\nabla \varphi$ is contained in $\{y\in\Rset^N: R\le |y-x|_*\le 2R\}$.

Let $L$ be the radial Laplacian operator in dimension $N-1$, that is,
$L(\phi)(r)\decl\phi''(r)+\frac{N-2}{r}\phi(r)$,
and set $t\decl\frac{|x_N-y_N|}{R},\ \ s\decl\frac{|x'-y'|}{R}$.
By computation we have,
$$\Delta\varphi(y)=\frac{1}{R^2}\left[\phi''(t)\phi(s)
+\phi(t)L(\phi)(s)\right],$$
hence,
$$|\Delta\varphi(y)|\le \frac{M}{R^2}.$$
Now, the supports of $\Delta\varphi$ and of $\nabla \varphi$ are such that,
$supp(\Delta\varphi) \subset supp(|\nabla \varphi|)\subset \{y:R\le |x-y|_*\le 2R\}.$
In what follows we denote by $A^*_R(x)$ the annulus
$A^*_R(x)\decl B^*_{2R}(x)\setminus B^*_R(x)\subset \CHS$.
See Figure \ref{fig1} at page \pageref{fig1}.

Therefore we have, 
$$  |I_2|\le \frac{M}{R^2} \int_{A^*_R(x)} |u|G^x_\epsilon(y) dy.$$
Observing that for $R\le\abss{x-y}\le 2R$ we have, 
$%2R^2\ge\epsilon^2+|\hat x-y|^2=\xg\ge\xp=
\epsilon^2+|x-y|^2>R^2,$ from 1. of Proposition~\ref{prop1}
we deduce 
$  |G^x_\epsilon|\le C'_n\frac{x_Ny_N}{R^{N}}
%\le  c\frac{x_Ny_N}{R^N},
$, hence,
\begin{equation}\label{estI2}
  |I_2|\le  \frac{c x_N}{R^{N+2}} \int_{A^*_R(x)} y_N |u(y)| dy.
\end{equation}

Next, we estimate $|I_3|$. We begin by computing
$|\nabla \varphi\cdot\nabla G^x_\epsilon|$ in $A^*_R(x)$ with
$R>x_N$.
Since
$$
\begin{aligned}
|I_3|\le& 2 \int_{A^*_R(x)}
|u(y)| |\nabla G^x_\epsilon\cdot\nabla\varphi|\\
\le& 
 2 \int_{A^*_R(x)}|u(y)| \abs{\sum_{j=1}^{N-1}\partial_j\varphi(y)\partial_j G_\epsilon^x(y)} +  
 2 \int_{A^*_R(x)}|u(y)| \abs{\partial_N\varphi(y)\partial_N G_\epsilon^x(y)}   
=:I_{31}+I_{32},
\end{aligned}
 $$
we get
$$
\begin{aligned}
 \abs{\sum_{j=1}^{N-1}\partial_j\varphi(y)\partial_j G_\epsilon^x(y)}=
\abs{\phi(\frac{\abs{x_N-y_N}}{R})\phi'(\frac{\abs{x'-y'}}{R})
\frac{1}{R} C'_N \abs{x'-y'}\left(\xg^{-N/2}-\xp^{-N/2} \right)}\\
\le\! c_1 \frac{\abs{x'-y'}}{R}\frac{\xg^{N/2}- \xp^{N/2}}{\xg^{N/2}\xp^{N/2}}
\le\! c_1 \frac{\abs{x'-y'}}{R}\frac{N}{2}(\xg-\xp)\xg^{N/2-1}\frac{1}{\xg^{N/2}\xp^{N/2}}
=\! c_1 \frac{\abs{x'-y'}}{R}\frac{N}{2} \frac{4x_Ny_N}{\xg\xp^{N/2}}.
\end{aligned}$$
Here, we have used the convexity inequality \eqref{conv} with $\alpha=N/2$, and \eqref{a1-a2}.
Therefore, since for $y\in A^*_R(x),$ we know that
$\xg\ge \xp\ge R^2$, we infer
$$\abs{\sum_{j=1}^{N-1}\partial_j\varphi(y)\partial_j G_\epsilon^x(y)}
\le c \frac{y_N}{R^{N+2}},$$
while, 
$$\abs{\partial_N\varphi(y)\partial_N G_\epsilon^x(y)} =  
\abs{\frac 1R\phi'(\frac{\abs{x_N-y_N}}{R})\phi(\frac{\abs{x'-y'}}{R})
 C'_N \left(\frac{y_N+x_N}{\xg^{N/2}}-\frac{y_N-x_N}{\xp^{N/2}} \right)}.
$$
Now, we observe that the above term can be non zero only for $R\le|x_N-y_N|\le 2R$.
Next by choosing $R\ge 2|x_N|$, it follows that $|y_N|\ge |x_N-y_N| -|x_N|\ge R/2,$
and $|y_N|\le|x_N-y_N| +|x_N|\le 5R/2$.
Hence for $\xg\ge\xp\ge R^2$, we obtain
$$\begin{aligned}
\abs{\partial_N\varphi(y)\partial_N G_\epsilon^x(y)}
\le c \left[ \frac{|y_N|}{R}  \left(\frac{1}{\xp^{N/2}}-\frac{1}{\xg^{N/2}}\right)+ \frac{|x_N|}{R}  \left(\frac{1}{\xp^{N/2}}+\frac{1}{\xg^{N/2}}\right) \right]\\
\le c\left[\frac 52 \frac{\xg^{N/2}-\xp^{N/2}}{\xg^{N/2}\xp^{N/2}}   
+ \frac{|x_N|}{R}  \left(\frac{1}{\xp^{N/2}}+\frac{1}{\xg^{N/2}}\right)\frac{2|y_N|}{R}\right]
\\ \le c\left[\frac 52 \frac{N}{2}(\xg-\xp)\xg^{N/2-1}\frac{1}{\xg^{N/2}\xp^{N/2}} + 2 x_N \frac{2}{R^N} \frac{y_N}{R^2}
\right]\\
 = c\left[\frac{5N}{4}\frac{4x_Ny_N}{\xg\xp^{N/2}}+ 2 x_N \frac{2}{R^N} \frac{y_N}{R^2} \right]
\le c x_N [5N+ 4] \frac{y_N}{R^{N+2}}.
\end{aligned}
$$

Therefore, we conclude that for $R>2x_N,$ we have
\begin{equation}\label{estI3}
  |I_3|\le  \frac{cx_N}{R^{N+2}} \int_{A^*_R(x)} y_N |u(y)| dy.
\end{equation}
It is important to note that that the constants $c$ in \eqref{estI2} and \eqref{estI3}
depends only on $N$ and on the cut-off function $\phi$.

To complete the proof
of Theorem \ref{teo:repiff}.\ref{mainA}, we first need 
 to take the limit as $\epsilon\to 0$  and then let  $R\to\infty$ in \eqref{repepsilon}.

Considering the right hand side of \eqref{repepsilon}, using the estimates above and  the assumption \rhz\!\!, it follows that
$$-\int_{\HS}u(y) \Delta( G^x_\epsilon \varphi)(y) dy\  \xrightarrow{\ \epsilon\to 0\ }
u(x)+{I_2}|_{\epsilon=0}+I_3|_{\epsilon=0}\xrightarrow{\ R\to \infty\ }u(x).$$

On the left  hand side of \eqref{repepsilon}, the limits
$$\int_{\DHS} K^x_\epsilon(y')\varphi(y',0)  d\nu(y')
 \xrightarrow{\ \epsilon\to 0\ }
\int_{\DHS} K^x(y')\varphi(y',0)  d\nu(y') 
\xrightarrow{\ R\to \infty\ }
\int_{\DHS} K^x(y')  d\nu(y'), 
$$
follow from the Beppo Levi monotone convergence theorem. From \eqref{repepsilon}
we deduce that these limits are finite.

The same reasoning applies to the first integral in \eqref{repepsilon},
yielding
$$ \int_{\HS}G^x_\epsilon(y)\varphi(y) d\mu(y) \xrightarrow{\ \epsilon\to 0\ }
 \int_{\HS}G^x(y)\varphi(y) d\mu(y)\xrightarrow{\ R\to \infty\ }
  \int_{\HS}G^x(y) d\mu(y).
$$
This completes the proof of \eqref{disrepmeas+}.

If $u$ solves \eqref{eqh}, then \eqref{repepsilon} holds with equality, and the proof follows by applying the same arguments.
\end{proof}

\subsection{Proof of Theorem \ref{repmeasnosign}}\label{sec:proof2}
We proceed using the same approach as in of the proof of Theorem \ref{teo:repiff}.\ref{mainA}.
The key difference in this proof lies in the conclusion, specifically in the computation of the limits of the following integrals:

\begin{equation}\label{intmeasures}
\int_{\DHS} K^x_\epsilon(y')\varphi(y',0)  d\nu(y') ,\quad and \quad \int_{\HS}G^x_\epsilon(y)\varphi(y) d\mu(y).
\end{equation}
Since $\nu=\nu^+-\nu^-$, the first integral can be expressed as
$$\int_{\DHS} K^x_\epsilon(y')\varphi(y',0)  d\nu(y') =\int_{\DHS} K^x_\epsilon(y')\varphi(y',0)  d\nu^+(y') -\int_{\DHS} K^x_\epsilon(y')\varphi(y',0)  d\nu^-(y') .$$
Each of these integrals can be analyzed separately.
% as we made in the end of proof of Theorem \ref{teo:repiff}.\ref{mainA}, 
By the Beppo Levi monotone convergence theorem, we have
$$\int_{\DHS}\! K^x_\epsilon(y')\varphi(y',0)  d\nu^\pm(y') 
 \xrightarrow{\ \epsilon\to 0\ }\!
\int_{\DHS}\! K^x(y')\varphi(y',0)  d\nu^\pm(y') 
\xrightarrow{\ R\to \infty\ }\!
\int_{\DHS}\! K^x(y')  d\nu^\pm(y').
$$
Given the assumption \eqref{hypneg}, it follows that
$$\int_{\DHS} K^x  d\nu(y') =\int_{\DHS} K^x d\nu^+(y') -\int_{\DHS} K^x d\nu^-(y')$$
is well-defined, and from \eqref{repepsilon}, we deduce that it is finite.

The second integral in \eqref{intmeasures} can be handled similarly.

\subsection{Proof of Theorem \ref{teo:repiff}.\ref{teo:inverse}}\label{proof1B}
Let $u$ be defined by \eqref{eqrepmeas+}.

To prove that $u\in L^1_{loc}(\CHS)$, we start by defining
\begin{equation}
    u_1(x)\decl \int _{\DHS} K^x(y') d\nu(y'),\quad   u_2(x)\decl \int_{\HS} G^x(y) d\mu(y),\qquad for\ x\in\HS.
 \end{equation}
It suffices to show that
\begin{equation}\label{ujfinito}
\int_{B^*_R} u_j(x) dx<\infty\quad for\ any\ R>0,\end{equation}
for $j=1,2$.

For $y'\in\Rset^{N-1},$ set $\displaystyle f_1(y')\decl \int_{B^*_R}K^x(y') dx$.
We then have
$$ \begin{aligned}
\int_{B^*_R} u_1(x)dx&= \int_{B^*_R}  \int _{\DHS} K^x(y') d\nu(y') dx=
 \int _{\DHS}\left(\int_{B^*_R}K^x(y') dx\right) d\nu(y')\\
 &=  \int _{\DHS}f_1(y') d\nu(y')=
   \int _{|y'|\le S}f_1(y') d\nu(y')+ \int _{|y'|> S}f_1(y') d\nu(y'),
   \end{aligned}
$$
where $S>0$ is chosen later.
Since $u_1$ is finite a.e., $K^x(y')\approx |y'|^{-N}$ as $|y'|\to \infty$,
and $\nu$ is finite on bounded sets, we deduce that  
\begin{equation}\label{nufinite}
  \int _{\DHS} \frac{1}{1+|y'|^N} d\nu(y')<\infty.\end{equation}
Next, we claim that $f_1(y')\approx |y'|^{-N}$  as $|y'|\to \infty$. Indeed,
applying the Lebesgue dominated convergence theorem 
and Proposition \ref{prop1}(3), we obtain
$$ \lim_{|y'|\to \infty}f_1(y')|y'|^{N}=
\lim_{|y'|\to \infty} \int_{B^*_R}K^x(y')|y'|^{N}  dx=\int_{B^*_R} C'_N x_N dx<\infty.$$
This fact, together with \eqref{nufinite}, yields,  
$\int _{|y'|> S}f_1(y') d\nu(y')<\infty$ for $S$ large.

To prove that the remaining integral is also finite, it is sufficient to show that $f_1$ is locally bounded.
Consider the expression
$$f_1(y')=C'_N\int_{|x'|<R}\int_0^R\frac{x_N}{(x_N^2+|x'-y'|^2)^{N/2}}dx_N dx'=
C'_N\int_{|x'|<R}\gamma(|x'-y'|)dx',$$
where we define $$\gamma(|x'-y'|)\decl \int_0^R\frac{x_N}{(x_N^2+|x'-y'|^2)^{N/2}}dx_N.$$
Thus, we have
$f_1=C'_N\chi_{B_R}\ast \gamma$, where the convolution is taken in $\Rset^{N-1}.$  Since $\chi_{B_R}\in L^\infty(\Rset^{N-1})$ and has bounded support,
if $\gamma\in L^1_{loc}(\Rset^{N-1})$,  then the claim $f_1\in  L_{loc}^\infty(\Rset^{N-1})$
follows by standard arguments.
Explicitly, in the case $N=2,$ we have
$$\gamma(|x'-y'|)=\int_0^R\frac{x_N}{(x_N^2+|x'-y'|^2)}dx_N=\frac12\ln
\frac{R^2+|x'-y'|^2}{|x'-y'|^2}.$$
For  $N\ge 3$, we have
$$\begin{aligned}
 \gamma(|x'-y'|)&=\int_0^R\frac{x_N}{(x_N^2+|x'-y'|^2)^{N/2}}dx_N\\
 &= \frac{1}{N-2}\left(\frac{1}{|x'-y'|^{N-2}}-\frac{1}{(R^2+|x'-y'|^{2})^{(N-2)/2}}\right).\end{aligned}
$$
This completes the proof of \eqref{ujfinito} for $j=1$.

The proof of \eqref{ujfinito} for $j=2$ follows the same steps as before.
For $y\in\HS,$ let $\displaystyle f_2(y)\decl \int_{B^*_R}G^x(y) dx$.
We have
$$ \begin{aligned}
\int_{B^*_R} u_2(x)dx&= \int_{B^*_R}  \int _{\HS} G^x(y) d\mu(y) dx=
 \int _{\HS}\left(\int_{B^*_R}G^x(y) dx\right) d\mu(y)\\
 &=  \int _{\HS}f_2(y) d\mu(y)=
   \int _{|y|\le S}f_2(y) d\mu(y)+ \int _{|y|> S}f_2(y) d\mu(y),
   \end{aligned}
$$
where $S>0$ will be chosen later.

Next, observe that $f_2(y)\approx y_N/|y|^N$ as $|y|\to \infty$. Indeed,
since $G^x(y)\approx y_N|y|^{-N}$ as $|y|\to \infty$
(see 2. of Proposition \ref{prop1}), by the dominated convergence theorem, we obtain
$$ \lim_{|y|\to \infty}f_2(y)|y|^{N}/y_N=
\lim_{|y|\to \infty} \int_{B^*_R}G^x(y)|y|^{N}/y_N  dx=\int_{B^*_R} C'_N x_Ndx<\infty.$$
Since 
\begin{equation}\label{mufinite}
  \int _{\HS} \frac{y_N}{1+|y|^N} d\mu(y)<\infty,\end{equation}
we have $$\int _{|y|> S}f_2(y) d\mu(y) <\infty\quad for\ S\ large.$$
Moreover, 
$$f_2(y)=\int_{B^*_R}G^x(y) dx=\int_{B^*_R}\Gamma^x(y) dx  -\int_{B^*_R}\Gamma^{\hat x}(y) dx=\chi_{B^*_R}\ast \Gamma^y-\chi_{B^*_R}\ast \Gamma^x(\hat y),
$$
where $f_2$  is the  sum of two functions which are convolutions of $\chi_{B^*_R}$
with locally integrable kernels. Hence, $f_2$ is locally bounded.
This concludes the proof of \eqref{ujfinito} for $j=2$.

\medskip

To prove that $u$ is a weak solution of \eqref{eqh},
 we choose $\varphi \in \D$
 % (the space of smooth, compactly supported function) . 
 and we aim to show that \eqref{eq:defsol} holds with equality sign.
 %, where $f$ and $g$ are functions defined in the context of the problem, and $u$ is the function we want to prove is a weak solution.
This requires us to compute the following integrals
$$\int_{\HS,x} \left( (-\Delta_x \varphi)(x)\int_{\HS,y}G^x(y)d\mu(y)\right)dx
+ \int_{\HS,x} \left(  (-\Delta_x \varphi)(x)\int_{\DHS,y'} K^x(y')d\nu(y')\right)dx.$$

Using Fubini-Tonelli theorem, the symmetry $G^x(y)=G^y(x)$ and the fact that $G$ is the Green function of $-\Delta$ on $\HS$, we compute the first integral as follows:
$$\begin{aligned}
\int_{\HS,x} \int_{\HS,y} (-\Delta_x \varphi)(x)G^x(y)d\mu(y)dx&=
\int_{\HS,y} \left(\int_{\HS,x} (-\Delta_x \varphi)(x)G^y(x) dx\right)d\mu(y)\\
&=\int_{\HS,y} \varphi(y)d\mu(y). \end{aligned}
$$
Since this integral is finite for any $\varphi \in \D$, we can conclude that 
 $\mu\in {\cal M}_{loc}(\CHS,x_N)$ (see Remark \ref{rem:2}).
 
Similarly, for the second integral, we have: 
$$
\begin{aligned}
\int_{\HS,x}  \int_{\DHS,y'} (-\Delta_x \varphi)(x)K^x(y')d\nu(y')dx\qquad\qquad\qquad\qquad\qquad\qquad\\=
 \int_{\DHS,y'}\left(\int_{\HS,x}  (-\Delta_x \varphi)(x)\frac{\partial}{\partial y_N} G(x,y)_{|y_N=0} dx\right)d\nu(y')\quad\\
  =\int_{\DHS,y'}\frac{\partial}{\partial y_N}
  \left( \int_{\HS,x}  (-\Delta_x \varphi)(x)G(x,y)dx\right)_{|y_N=0}d\nu(y')\quad\\
 =  \int_{\DHS,y'}\frac{\partial}{\partial y_N}\left(\varphi(y) \right)_{|y_N=0}d\nu(y')
  = \int_{\DHS,y'}\varphi_N(y',0) d\nu(y').
 \end{aligned}
 $$
  Here $\phi_N(y',0)$ denotes the normal derivative of $\phi$ at $y_N=0.$  By summing the results of these integrals, we conclude that $u$ satisfies the weak formulation of the equation \eqref{eqh}, confirming that $u$ is indeed a weak solution.

\medskip

To complete the proof of Theorem \ref{teo:repiff}.\ref{teo:inverse},
we turn now to the transformation proposed by Huber \cite{H56}, following an idea of Weinstein.  
Define a transformation of the function $v(\xi) = H[u(x)]$ in the extended space $\Rset^{N+2},$ 
where $\xi=(\xi_1,\dots,\xi_{N-1},\xi_N,\xi_{N+1},\xi_{N+2})=(\xi',\bar\xi)\in\Rset^{N+2}$, with the notation $\xi'=(\xi_1,\dots,\xi_{N-1})\in \Rset^{N-1}$,  $\bar \xi=(\xi_N,\xi_{N+1},\xi_{N+2})\in \Rset^{3}$,
and set $$S\decl \{\xi\in \Rset^{N+2}, |\bar\xi|=0\}.$$
Let $v(\xi)=H[u(x)]$  be defined on $\Rset^{N+2}$ by
%\begin{equation}\label{transHuber}
%\begin{aligned}
%v(\xi)&\decl\frac{u(\xi_1,\xi_2,\dots,\xi_{n-1},(\xi_n^2+\xi_{n+1}^2+\xi_{n+2}^2)^{1/2})}{(\xi_n^2+\xi_{n+1}^2+\xi_{n+2}^2)^{1/2}},\quad  for\ \ \xi\not\in S,\\
%&and\ by\\
%v(\xi)&\decl\liminf_{\eta\to\xi, \eta\not\in S} v(\eta),\qquad\qquad\qquad\qquad%\qquad for  \ \ \xi\in S.
%\end{aligned}
%\end{equation}

\begin{equation}\label{transHuber}
v(\xi)\decl \left\{\begin{array}{ll}
\displaystyle \frac{u(\xi_1,\xi_2,\dots,\xi_{n-1},(\xi_n^2+\xi_{n+1}^2+\xi_{n+2}^2)^{1/2})}{(\xi_n^2+\xi_{n+1}^2+\xi_{n+2}^2)^{1/2}} &\qquad  for\ \ \xi\not\in S,\\ \\  \displaystyle
\liminf_{\eta\to\xi, \eta\not\in S} v(\eta) &\qquad for  \ \ \xi\in S.
\end{array}\right.
\end{equation}
\begin{lemma}[\cite{H56}]\label{lem:huber} If $u$ is superharmonic on $\HS$ and
  \begin{equation}\label{HuberBound}
   \liminf_{x\to y} u(x)\ge 0 \quad for\ all\ y\in\DHS,
  \end{equation}
  then $v$ is superharmonic on $\Rset^{N+2}$.
  
Conversely,  if $v$ is superharmonic on $\Rset^{N+2}$ and is symmetric with respect to
  the subspace $S$ \footnote{This means that it depends only on 
  $\xi_1,\xi_2,\dots,\xi_{N-1}$ and $(\xi_N^2+\xi_{N+1}^2+\xi_{N+2}^2)^{1/2}$.},
  then the function,
  $$u(x)=x_Nv(x_1,\dots,x_{N-1},x_N,0,0)$$
   is superharmonic on $\HS$ and satisfies
  \eqref{HuberBound}.
\end{lemma}
Notice that in the Lemma above, there is no assumption on the sign of $u$ (or, equivalently, of $v$).

\subsubsection*{Completion of the Proof of Theorem \ref{teo:repiff}.\ref{teo:inverse}}
We need to show that  \rhh holds.
By replacing $u(x)$ with $u(x)-hx_N$, i.e. assuming that $h=0$,  it follows that $u$ is a nonnegative and superharmonic weak solution of \eqref{eqh}.

Let $ v(\xi)=H[u]$.
Since $u$ is nonnegative, we deduce that $v$ is a nonnegative superharmonic function in $\Rset^{N+2}$.
Let $c\decl\inf v\ge 0$. 
It is known that (see \cite{CDM})
$$ \lim_R \frac{1}{R^{N+2}}\int_{B^{e,n+2}_R(\eta_0)} (v(\xi)-c)d\xi =0,\quad for\ a.e.\  \eta_0\in \Rset^{N+2},$$
where
$B^{e,n+2}_R(\eta_0)$ is the Euclidean ball in $\Rset^{N+2}$ of radius $R$
centered ad $\eta_0\in\Rset^{N+2}$.
In the integral above, the Euclidean ball can be replaced by the ball 
$B^{*,N+2}_R(\eta_0)$ associated with the norm
$$|\xi|_*=|(\xi_1,\dots,\xi_{N-1},\xi_N,\xi_{N+1},\xi_{N+2})|_*=
|(\xi',\bar\xi)|_*\decl
\max\{ |\xi'|,|\bar\xi|\}.
$$
Here, we have used a notation similar to the one we have used in $\HS$.
%introduced in  Appendix~\ref{sec2}.
Now, fix $x=(x',x_N)\in\HS$ and set $\eta_0\decl(\eta_0',\bar\eta_0)=(x',x_N,0,0)$,
so that
$\eta_0'=x'$, $\bar\eta_0=(x_N,0,0)$ and $|\bar\eta_0|=x_N$.
Therefore, we have
\begin{equation}\label{liminfv}
\lim_R \frac{1}{R^{N+2}}
\int_{B^{*,N+2}_R(\eta_0)} (v(\xi)-c)d\xi =\lim_R \frac{1}{R^{N+2}}
\int_{|\xi'-\eta'_0|<R} \int_{ |\bar\xi-\bar\eta_0|<R} (v(\xi)-c)d\bar\xi d\xi'=0.\end{equation}
Let us estimate the last integral for large $R$.
Let $R\ge 2\abs{\bar\eta_0}$ and set $\delta\decl \abs{\bar\eta_0}/R$,
$\tau\decl \sqrt{1-\delta^2}$. 
In $\Rset^3$, with this choice, we have the inclusion
$$\left\{\bar\xi:\abs{\bar\xi}<\tau R,\xi_N>0\right\}\subset B^{e,3}_R(\bar\eta_0)=
  \{\bar\xi:\abs{\bar\xi-\bar\eta_0}<R\}.$$
Recalling that $v(\xi)=v(\xi',|\bar\xi|)$, we obtain:
 $$\begin{aligned}
   \int_{ |\bar\xi-\bar\eta_0|<R} |v(\xi)-c|d\bar\xi \ge 
     \int_{ \xi_N>0\atop |\bar\xi|<\tau R} |v(\xi',|\bar\xi|)-c|d\bar\xi =
     \frac{\sigma_3}{2}\int_0^{\tau R} |v(\xi',|r)-c|r^2d r \qquad\qquad\\
    = \frac{\sigma_3}{2}\int_0^{\tau R} |u(\xi',r)-cr|r d r 
    \ge    \frac{\sigma_3}{2}\int_{y_N>0\atop|y_N-x_N|<\gamma R} |u(\xi',y_N)-cy_N|y_N d y_N 
 \end{aligned}$$
where $\gamma>0$ is such that $\{y_N:|y_N-x_N|<\gamma R,y_N>0\}\subset [0,\tau R]$. This is possible whenever $\gamma R+x_N\le \tau R,$ which holds for
$0<\gamma\le \tau-\delta=\sqrt{1-\delta^2}-\delta$.

Therefore,  for any $0<\gamma\le \sqrt{1-\frac{x_N^2}{R^2}}-\frac{x_N}{R}$, or equivalently
for any $0<\gamma<1$ and $R>2\max\{1,(\sqrt{1-\gamma^2}-\gamma)^{-1}\} x_N,$ \footnote{For instance, choosing $\gamma=1/2$, \eqref{intliftedcorr} holds for any  $R>\frac52x_N$.}, we have
\begin{equation}\label{intliftedcorr}
\begin{aligned}\int_{B^{*,N+2}_R(\eta_0)} |v(\xi)-c|d\xi &=\int_{|\xi'-\eta'_0|<R} \int_{ |\bar\xi-\bar\eta_0|<R} |v(\xi)-c|d\bar\xi d\xi' \\
&\ge \frac{\sigma_3}{2}
\int_{|\xi'-\eta'_0|<\gamma R} d\xi' \int_{y_N>0\atop|y_N-x_N|<\gamma R} |u(\xi',y_N)-cy_N|y_N d y_N \\
&=\frac{\sigma_3}{2}\int_{y\in B^*_{\gamma R}(x)\cap \HS} |u(y)-cy_N|y_N dy.\end{aligned}
\end{equation}

Relation \eqref{intliftedcorr}, combined with \eqref{liminfv} and the fact that $v\ge c,$ implies by  rescaling  the parameter $R$
that $u$ satisfies \rhh with a constant $c\ge 0$. We aim to prove that $c=0$.
Since $u$ is a weak solution of \eqref{eqh} and satisfies \rhh with a constant $c\ge 0$, 
$u$ can be written as: 
 \begin{equation*}
      u(x)= cx_N+\int _{\DHS} K^x(y') d\nu(y') +  \int_{\HS} G^x(y) d\mu(y).
   \end{equation*}
However, given that $u$ is defined by \eqref{eqrepmeas+} with $h=0$,
it follows that $c=0$.
Finally, we observe that \eqref{intliftedcorr}, along with \eqref{liminfv} and the fact that  and $v\ge c$, confirm the statement \ref{teo:repiff}.\ref{teo:inverse}.(d).

{Proof of \ref{teo:repiff}.\ref{teo:inverse}.(a).}
From the above argument, it is evident that
the infimum of $v$ is the constant $h$ appearing in the representation of  $u$.

Therefore, by definition of $v,$ we have:
$$h=\inf_{\Rset^{N+2}} v =\inf_{\Rset^{N+2}} \frac{u(\xi',|\bar \xi|)}{|\bar \xi|}.
$$

 {Proof of \ref{teo:repiff}.\ref{teo:inverse}.(b).}
Let $\tilde\Omega\subset \Rset^{N+2}$ be
defined as  $\tilde{\Omega}\decl\{\xi\in\Rset^{N+2}:(\xi',|\bar \xi|) \in\Omega\}$.
If $h=\inf_{\Omega} \frac{u(x)}{x_N}$, then  $h=\inf_{\tilde \Omega} v$.
Since $v$ is superharmonic on the whole $\Rset^{N+2}$ and $\tilde \Omega$ is bounded, it follows that $v\equiv h,$ thus  the claim is proved.

 {Proof of \ref{teo:repiff}.\ref{teo:inverse}.(c).}
To establish the estimate \eqref{est:ubelow}, we begin by noting that
at least one of the measures $\mu$ or $\nu$ must be non trivial.
Without loss of generality we may assume that $h=0$.
We first consider the case where $\mu\not\equiv 0$ and let $R_0$ be such that
$\mu(B^*_{R_0})>0$.
Since $G^x(y)=G^y(x)$, from part 2 of Proposition \ref{prop1}, it follows that, $G^x(y)\ge c\frac{x_N}{|x|^N}$ 
for large $|x|$ and
for any $y\in B^*_{R_0}$. 
Therefore, 
$$ u(x)\ge \int_{\HS}G^x(y)d\mu(y)\ge  \int_{B^*_{R_0}}G^x(y)d\mu(y)
\ge \int_{B^*_{R_0}} c\frac{x_N}{|x|^N} d\mu(y)\ge  c \mu(B^*_{R_0})\frac{x_N}{|x|^N}.$$

Now consider the case $\nu\not\equiv 0$ and let $R_0$ be such that
$\nu(B'_{R_0})>0$.
For $y'\in B'_{R_0}$, we have 
$$|x'-y'|\le |x'|+R_0\le \sqrt2\sqrt{ |x'|^2+R_0^2}.$$ Hence,
$$x_N^2+|x'-y'|^2\le 2 (x_N^2+|x'|^2+R_0^2),$$ which implies
$$ K^x(y')\ge c \frac{x_N}{(x_N^2+|x'|^2+R_0^2)^{N/2}},\quad for\ |y'|<R_0.$$
Therefore, 
$$ u(x)\ge \int_{\DHS} K^x(y')d\nu(y')\ge  \int_{B'_{R_0}} c \frac{x_N}{(x_N^2+|x'|^2+R_0^2)^{N/2}} d\nu(y')= c \nu(B'_{R_0})\frac{x_N}{(|x|^2+R_0^2)^{N/2}}.$$
In summary, the claim is proved for $|x|$ large, say for $|x|>R_1$.

Since  \ref{teo:repiff}.\ref{teo:inverse}.(a) and  \ref{teo:repiff}.\ref{teo:inverse}.(b) hold, the claim for any $x$ follows by observing that, 
$\inf_{B_{R_1}} \frac{u(x)}{x_N}>0=h$. 

\subsection{Proof of Theorem \ref{teo:potentialink}}\label{secproofpotlink}

\begin{proof}
Let $w\decl u-v$. Since $-\Delta w\ge 0$ in the distributional sense, there exists a positive Radon measure $\mu$ such that $-\Delta w =\mu$ and
$$\liminf_{x\to (y',0)} w(x)\ge 0,\quad for\ \ y'\in\DHS.$$

Using the same notations as on page \pageref{transHuber},
set $z\decl H[w]$ as in \eqref{transHuber}. 
From Lemma~\ref{lem:huber}, it follows that $z$ is superharmonic in $\Rset^{N+2}$.
Fix $x=(x',x_N)\in \HS$  and let
$\eta_0= (\eta_0',\bar\eta_0)\decl(x',x_N,0,0)$.
Arguing as in \eqref{intliftedcorr} for $R>\frac52 x_N$, we have 
$$
\int_{B^{*,N+2}_{2R}(\eta_0)} |z(\xi)|d\xi\ge\frac{\sigma_3}2\int_{B^*_R(x)\cap\HS} |w(y)| y_N dy.$$
Now, since $z$ is superharmonic,  we know that $z\in L^1_{loc}(\Rset^{N+2})$.
Hence we deduce that the above integrals are finite and
$x_N|w|\in  L^1_{loc}(\CHS)$. This completes the proof of the claim.

Next we prove the equivalence between \eqref{umv} and \eqref{u-vRing}, i.e., the equivalence between $w\ge 0$ and the fact that $w$ satisfies \rhh with $h\ge 0$.

Step 1 : Assume $w\ge 0$. Since $w$ is a superharmonic function, by 
Theorem~\ref{th:RepHSClassic}, it follows $w$ can be represented by \eqref{eqrepmeas+}.
An application of Theorem \ref{teo:repiff}.\ref{teo:inverse} completes the proof of the claim, implying $w\in L^1_{loc}(\CHS)$.

Step 2: Assume $w$ satisfies \rhh with $h\ge 0$. 
%By replacing $w$ with $w-hx_N$, 
Without loss of generality we may suppose that $w$ satisfies \rhz\!\!. 
Notice that if $z\decl H[w]$ is nonnegative, then $w\ge0$. 
To prove that $z\ge 0$, we apply Theorem \ref{th:RepClassic} to the function $z$ with $l=0$ in $\Rset^{N+2}$.
%{\color{red}---DISCORSO CHE (questo e' vero anche per gli anelli nel semispazio, mettiamo in un remark e qui lo richiamiamo)----}
Let $1>\tau>1/\sqrt{2}$ be fixed. Consider $x=(x',x_N)\in \HS$  and let
$\eta_0= (\eta_0',\bar\eta_0)\decl(x',x_N,0,0)$. Our claim is that: 
\begin{equation*}
\liminf\frac{1}{R^{N+2}}\int_{\frac{\sqrt2 R}\tau<|\eta_0-\xi|<2R} |z(\xi)|\,d\xi=0.
\footnote{Obviously, in condition \ringentire the annulus $B_{2R}\setminus B_R$ can be replaced by the annulus $B_{\gamma R}\setminus B_R$ with any $\gamma>1$.}
\end{equation*}
Noticing that $B^{*,N+2}_{2R}(\eta_0)\supset B^{e,N+2}_{2R}(\eta_0)$
and $B^{*,N+2}_{R}(\eta_0)\subset B^{e,N+2}_{\sqrt2R}(\eta_0),$
 we have $$\int_{\frac{\sqrt2 R}\tau<|\eta_0-\xi|<2R} |z(\xi)|\,d\xi\le 
  \int_{B^{*,N+2}_{2R}(\eta_0)\setminus B^{*,N+2}_{R/\tau}(\eta_0)} |z(\xi)|\,d\xi.$$
Therefore, it suffices to prove that for $R>\frac{2\tau}{1-\tau} x_N,$ there holds:
\begin{equation}\label{dis:tec1}
\int_{B^{*,N+2}_{2R}(\eta_0)\setminus B^{*,N+2}_{R/\tau}(\eta_0)} |z(\xi)|\,d\xi \le 2\sigma_3
\int_{y_N>0\atop R<\abss{x-y}<2R}  |u(y)|y_N \ dy.\end{equation}
The claim will follow from the hypothesis \rhz\!\!.

Step 3: Proof of \eqref{dis:tec1}. We  split the integration domain as
$B^{*,N+2}_{2R}(\eta_0)\setminus B^{*,N+2}_{R/\tau}(\eta_0)=A\cup B$
where:
$A\decl\{|\bar \xi-\bar\eta_0|<R/\tau, R/\tau\le |\xi'-\eta'_0|<2R\}$
and $B\decl\{R/\tau\le|\bar \xi-\bar\eta_0|<2R, |\xi'-\eta'_0|<2R\}$.
Since 
\begin{equation}\label{intA}
\int_A |z(\xi)|d\xi = \int_{R/\tau\le |\xi'-\eta'_0|<2R}d\xi' 
\int_{|\bar \xi-\bar\eta_0|<R/\tau} |z(\xi',\bar\xi)|d\bar\xi,
\end{equation}
and considering that $z(\xi',\bar\xi)=z(\xi',|\bar\xi|)$, we have:
$$\begin{aligned}
\int_{|\bar \xi-\bar\eta_0|<R/\tau} |z(\xi',\bar\xi)|d\bar\xi\le 
\int_{|\bar \xi|<R/\tau+|\eta_0|} |z(\xi',\bar\xi)|d\bar\xi=
\sigma_3 \int_0^{R/\tau+|\eta_0|} |z(\xi',r)|r^2dr\qquad\qquad\\
=\sigma_3 \int_{y_N>0\atop |y_N-x_N|<R/\tau}|u(\xi',y_N)|y_N dy_N
\le \sigma_3 \int_{y_N>0\atop |y_N-x_N|<2R}|u(\xi',y_N)|y_N dy_N,
\end{aligned}
$$
where in the last identity, we have used the fact that $R>\tau x_N$.
Therefore, from \eqref{intA}, we can conclude that:
$$\int_A |z(\xi)|d\xi \le \sigma_3\!\! \int_{R\le |\xi'-\eta'_0|<2R}d\xi'\!\! \int_{y_N>0 \atop |y_N-x_N|<2R}\!\!|u(\xi',y_N)|y_N dy_N\le 
\sigma_3\!\! \int_{y_N>0\atop R<\abss{x-y}<2R}\!\!  |u(y)|y_N \, dy.
$$
 Estimate for region $B$: Similarly, for region $B$, we have:
\begin{equation}\label{intB}
\int_B |z(\xi)|d\xi = \int_{ |\xi'-\eta'_0|<2R}d\xi' 
\int_{R/\tau\le |\bar \xi-\bar\eta_0|<2R} |z(\xi',\bar\xi)|d\bar\xi.
\end{equation}
Since
$\{R/\tau\le |\bar \xi-\bar\eta_0|<2R\}\subset \{R/\tau -x_N< |\bar \xi|<2R+x_N\}$, we obtain:
$$\begin{aligned}\int_{R/\tau\le |\bar \xi-\bar\eta_0|<2R} |z(\xi',\bar\xi)|d\bar\xi
\le \int_{R/\tau -x_N< |\bar \xi|<2R+x_N} |z(\xi',\bar\xi)|d\bar\xi\qquad\qquad\qquad\qquad\qquad\qquad\\
=\sigma_3\int_{R/\tau -x_N}^{2R+x_N}|z(\xi',r)|r^2 dr=
\sigma_3\int_{R/\tau -2x_N<y_N-x_N<2R}|u(\xi',y_N)|y_N dy_N\\
\le \sigma_3\int_{R<y_N-x_N<2R}|u(\xi',y_N)|y_N dy_N,\qquad\qquad\qquad\qquad\qquad\qquad\qquad
\end{aligned}$$
where we used the fact that $R>2\tau x_N$
and $R>\frac{2\tau}{1-\tau} x_N$.
Therefore, from \eqref{intB}, we can conclude that:
$$\int_B |z(\xi)|d\xi \le \sigma_3\! \int_{ |\xi'-\eta'_0|<2R}d\xi' \!\!\int_{R<y_N-x_N<2R}\! |u(\xi',y_N)|y_N dy_N\le 
\sigma_3\! \int_{y_N>0\atop R<\abss{x-y}<2R} \!\! |u(y)|y_N \, dy.
$$
This completes the proof.
\end{proof}

\section{Another integral representation formula}\label{sec:ringreen}
In this section we prove some representation formulae concerning regular superharmonic functions in the half-space $\HS$. These results can be deduced from some theorems contained in [11].
In what follows  $\mu$ stands for a continuous function defined on $\HS$.

The possibility to represent a superharmonic function $u$ in the whole space 
$\Rset^N$ or on the half-space $\HS$ with an integral formula, is linked to the asymptotic behavior of
some weighted integrals of the function $u$ on some rings, see Theorem \ref{th:RepClassic} and \ref{teo:repiff} respectively. 

Let $x\in\HS$ and $r>0$. We set 
$$\Omega_r(x):=\left\{y\in\HS \, |\, G^x(y)>\frac{1}{r} \right\}\cup \{x\}.$$

Throughout this section we shall call the set
$\Omega_{2R}(x)\setminus \Omega_R(x)$ the standard annulus.
 See Figure  \ref{fig1}.
\begin{figure} 
\captionsetup{width=0.7\textwidth}
\centering
%\ \ \ \ \ \ \ \ 
\includegraphics[scale=.7]{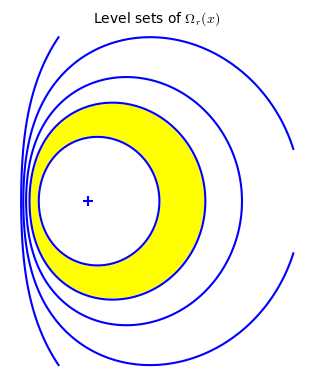}
\ \ \ \ \ \ \ \ \includegraphics[scale=.7]{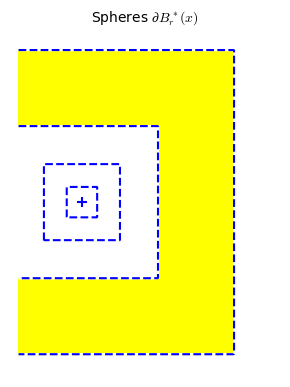}
\caption{On left: In blue the level sets  $\partial\Omega_r(x)$ with $x=(0,5)$ and $r=1,2,3,4,5$. The yellow filled region is the "annulus" $\Omega_2(x)\setminus\Omega_1(x)$.\newline
On right: In dotted line the surface of the ball $B^*_r(0,5)$ 
(in the $\abss \cdot$ norm) for $r=1,3,6,12$ 
The yellow filled region is the "annulus" $A^*_{6}(x)=B^*_{12}(x)\setminus B^*_{6}(x)$.}\label{fig1}
\end{figure}

Notice also that the rings of condition 
$({\cal R})$ are modeled in a similar way modulo a rescaling. Indeed
the integration domain appearing in $({\cal R})$ is given by $B_{2R}(x)\setminus B_R(x)$. Notice that
$$ B_r(x)=\left\{y\in\HS \, |\, \Gamma^x(y)>\frac{C_N}{r^{N-2}} \right\}\cup \{x\}.$$

\begin{theorem}\label{teo:repold} Let $u\in \C^2(\HS) $ 
  be such that $-\Delta u=:\mu\ge 0$. 
 
\begin{enumerate}
  \item\label{pointA}
  Let $x\in\HS$ and assume that
     \begin{equation}\label{eq:ringbist}
     l_x\decl \frac{1}{\ln2}\liminf_{R\to +\infty}
      \int_{\Omega_{2R}(x)\setminus\Omega_{R}(x)}
        \frac{|\nabla G^x(y)|^2} {G^x(y)}u(y) dy\in\Rset,
     \end{equation}
  then
  \begin{equation}\label{eq:repx}
   u(x)=l_x+\int_{\HS}G^x(y)\mu(y) dy.  
  \end{equation}
  
Therefore assuming that  (\ref{eq:ringbist}) holds for any $x\in\HS$, and setting
  $l(x)\decl l_x$, we deduce that, 
%  Assume now that for any $x\in\HS$ assumption  (\ref{eq:ringbist}) holds and set  $l(x)\decl l_x$. Then
  \begin{enumerate}
  \item\label{pointA1}  $\inf u = \inf l$ (finite or infinite),
  and  the following alternative holds,
  \begin{equation}    \label{eq:alt} either\ \ 
     u(x) > l(x), \forall x\in\HS,  \qquad  or \qquad \
        u\equiv l\ \ and \ \ \mu\equiv 0.  
   \end{equation}  
  \item\label{pointA2}  $l$ is  harmonic in $\HS$.
  \item\label{pointA4} If $l_x$ 
  does not depends on $x$, i.e.
  $l_x=l\in\Rset$, then  for any  $x\in\HS$ 
  \begin{equation}\label{eq:forrep}
    u(x)=l+\int_{\HS}G^x(y)\mu(y) dy.
  \end{equation}
  \end{enumerate}

 \item\label{pointB} If $u$ is bounded from below, then (\ref{eq:ringbist}) is fulfilled for any $x\in\HS.$ Hence
    the claims in \ref{pointA}. hold true.

 \item \label{pointC}
   If there exists a harmonic function $H\in \C^2(\HS)$ such that
   $$u(x)=H(x)+\int_{\HS}G^x(y)\mu(y) dy,$$
   then for any $x\in \HS$
   \begin{equation} \label{eq:ringh} H(x)=\frac1{\ln 2}\liminf_{R\to +\infty}
      \int_{\Omega_{2R}(x)\setminus\Omega_{R}(x)}
        \frac{|\nabla G^x(y)|^2} {G^x(y)}u(y)dy.
     \end{equation}
   \item\label{pointD} Let $x\in \HS$ and $l_x\in\Rset$. Then
   \begin{equation}\label{eq:ringbistVA}
	\liminf_{R\to +\infty}
      \int_{\Omega_{2R}(x)\setminus\Omega_{R}(x)}
        \frac{|\nabla G^x(y)|^2} {G^x(y)}|u(y)-l_x| dy=0,
     \end{equation}   
    if and only if
    \begin{equation} \label{eq:ringD}  \liminf_{R\to +\infty}
      R\int_{\Omega_{2R}(x)\setminus\Omega_{R}(x)}
		\frac{1}{|x-y|^{2N}}    |u(y)-l_x|=0.
     \end{equation}
    Moreover,  if one of these assumptions is satisfied, then
   \eqref{eq:ringbist} holds.
  \end{enumerate}
\end{theorem}
\begin{remark}
The existence of  a harmonic function $H$ in the hypothesis of point \ref{pointC}. above, 
is guaranteed by the Riesz decomposition theorem provided there exists a subharmonic minorant of $u$.
See also \cite[Theorem 4.4.1]{AG} and the comment preceding 
Theorem \ref{teo:distr=weak}.
\end{remark}

The following result is a consequence of Theorem 14 and Remark 17 of \cite{DM19}. 
Set
$$X\decl\{u\in \C^2(\CHS):u(x',0)=0,\ \ and\ \  \liminf_{x_N\to\infty} \frac{u(x)}{x_N}=0\}.$$
\begin{theorem}\label{teo:repgenold} The following statements hold.
  \begin{enumerate}
  \item Let  $u\in X$.
    Then $-\Delta u=:\mu\ge 0$ 
    and $\inf  u= 0$ if and only if
    $$u(x)=\int_{\HS} G^x(y)\mu(y) dy \qquad  \forall x\in\HS.$$
  \item  Let $u\in X$ be such that $-\Delta u\ge 0$.
    Then, $\inf u=0$    if and only if  
    \begin{equation}\label{eq:ringt}\liminf_{R\to +\infty}
      \int_{\Omega_{2R}(x)\setminus\Omega_{R}(x)}\frac{|\nabla G^x|^2}{G^x}
       u dy=0 \qquad  \forall x\in\HS,\end{equation}
     if and only if  
    \begin{equation}\lim_{R\to +\infty}
      \int_{\Omega_{2R}(x)\setminus\Omega_{R}(x)}\frac{|\nabla G^x|^2}{G^x}
      |u| dy=0 \qquad  \forall x\in\HS,\end{equation}
%  [resp. 
  if and only if 
    \begin{equation}\label{eq:ring0t}\liminf_{R\to +\infty}
      \int_{\partial\Omega_{R}(x)} |\nabla G^x|\ 
       u dH_{n-1}=0 \qquad \forall x\in\HS,\end{equation}
    if and only if  
    \begin{equation}\label{eq:ring0bist}\lim_{R\to +\infty}
      \int_{\partial\Omega_{R}(x)} |\nabla G^x|\ 
      |u-l| dH_{n-1}=0 \qquad \forall x\in\HS,\end{equation}
    if and only if  
    \begin{equation}
l=\lim_{R\to +\infty}
      \frac{1+\alpha}{R^{1+\alpha}} \int_{\Omega_{R}(x)}
      \frac{|\nabla G^x|^2}{({G^x})^{2+\alpha}}
      {u} dy\quad \alpha>-1\ \ \forall x\in\HS.
\end{equation}
  \item  Let $u\in X.$  Suppose that $u$ is bounded from below and $-\Delta u\ge 0.$ 
 Set
    \begin{equation}\label{eq:ringtalt}c_x:=\liminf_{R\to +\infty}\frac{1}{\ln 2}
      \int_{\Omega_{2R}(x)\setminus\Omega_{R}(x)}\frac{|\nabla G^x|^2}{G^x}\ u dy \qquad  \forall x\in\HS,\end{equation}
then $c_x$ does not depend on $x$ and  $c_x=0.$
  \end{enumerate}
\end{theorem}

\begin{proof}[Proof of Theorem \ref{teo:repold}]
   The theorem is a consequence of Remark 17 and Theorem 6 of \cite{DM19}.
   Indeed, the Green function on the half-space $\HS$ satisfies the assumptions 
   listed in \cite{DM19} H1.--H7. 
      
   The statements \ref{pointA}, \ref{pointA1},  \ref{pointA4},
   \ref{pointB} and \ref{pointC} are a direct consequence of
%     the corresponding statements in 
Theorem 6 in    \cite{DM19}.
   
The claim in \ref{pointA2} follows again from Theorem 6 in
   \cite{DM19} provided 
   $${w(x)\decl\int_{\HS}G^x(y)\mu(y) dy\in L^1_{loc}(\HS)}.$$ 
   This is easily seen to be true since $w$ is a superharmonic function.

 Proof of \ref{pointD}. First we observe that
the equivalence between \eqref{eq:ringbistVA} and \eqref{eq:ringD} is a consequence of the following argument.
For $x,y\in\HS$, $x\neq y$, 
from the estimates \eqref{est:GRabove2} and \eqref{est:GRbelow}
we have
$$ \left(\frac{C'_N}{2}\right)^2x_N^2\le \abs{\nabla G^x(y)}^2\abs{x-y}^{2N}
\le {C'_N}^2x_N^2[(1+N)^2+N^2].
$$
Let $x\in\HS$ be fixed, let $R>0$ and consider $\Omega_{2R}(x)\setminus\Omega_{R}(x)$. If $y\in \Omega_{2R}(x)\setminus\Omega_{R}(x)$
 we have
${2R}>\frac1{G^x(y)}>{R}$. 
Hence, for  $y\in \Omega_{2R}(x)\setminus\Omega_{R}(x)$ we deduce
$$\frac{ \abs{\nabla G^x(y)}^2}{G^x(y)}\approx \frac R{\abs{x-y}^{2N}}.$$
This completes the proof of  the equivalence  between \eqref{eq:ringbistVA} and \eqref{eq:ringD}.

Finally, the implication 
 \eqref{eq:ringbistVA} $\Rightarrow$ \eqref{eq:ringbist} follows from Lemma 26 in \cite{DM19}.
\end{proof}

\begin{proof}[Proof of Theorem \ref{teo:repgenold}] From Theorem 14 and Remark 17 of \cite{DM19} it follows that
the statements A., B. and C. are equivalent to the following classical Liouville property:

\textit{Let $u\in X$. If $u$ is harmonic and bounded from below, then $u$ is constant.}

\noindent   Indeed let $c>0$ be such that $u+c> 0$.
 By Theorem \ref{th:RepHSClassic} it follows  that 
$$u(x)+c=hx_N+\int _{\DHS} K^x(y') d\nu(y'),$$
where $\nu$ is a suitable positive Radon measure $\nu$ and $h\ge 0$.
Since $u\in X$ we have $\liminf_{x_N\to\infty}u(x)/x_N=0.$ 
Hence $h=0.$ In addition, by the fact that  $u(x',x_N)\to 0$ as $x_N\to 0$ it follows  that
$$ \int _{\DHS} K^x(y') d\nu(y')\to c\qquad as \ \ x_N\to 0,\ \ \forall x'\in \Rset^{N-1}.$$
In other words the harmonic function $w(x)\decl \int _{\DHS} K^x(y') d\nu(y')$
has boundary value identically equal to the constants $c.$ Therefore
$w\equiv c$ on $\HS,$ that is, $u\equiv 0$.
%Now since
%$$\int _{\Rset^{N-1}} \frac{C'_N x_N}{(x_N^2+|x'-y'|^2)^{N/2}} d y' =1 ,$$
%we deduce
%that$\nu\equiv c$, that is $u\equiv 0$.
\end{proof}

\renewcommand\appendixtocname{Appendix}
\renewcommand\appendixname{Appendix}
\begin{appendices}
\renewcommand\appendixname{Appendix}

\section{Estimates related to Poisson kernel and Green function}\label{sec2}
Here we collect a few results related to the fundamental solution 
of $\-\Delta$ on $\HS$, to the
Poisson kernel, and its Green function.

Let $R>0$ and set,
$$ B^*_R(x)\decl \left\{y\in\HS: \max\{|y_N-x_N|,|y'-x'|\}<R   \right\},$$
the "cylindrical" ball in $\HS$  centered at $x\in \HS$ and radius $R$ generated by the norm
$\abs\cdot_\ast$ defined in \eqref{def:norm}.
%\begin{equation}\label{def:norm}
%|x|_\ast=|(x',x_N)|_\ast\decl \max \{|x'|,|x_N|\},\end{equation}
%where $|\cdot|$ stands for the Euclidean norm.
%For $x\in\HS$ and $R>0$ we set
%$$ A^*_R(x)\decl B^*_{2R}(x)\setminus B^*_R(x).$$
%See Figure \ref{fig1}.

In this paper  $G$ and $K$ are defined as follows.
Let  $\Gamma^x(y)$ be the fundamental solution 
of $-\Delta$ at $x$, that is, 
$-\Delta_y \Gamma^x(y) = \delta_x$ where $\delta_x$ is the Dirac 
measure at $x$.
%Set, $C_N^{-1}\decl \sigma_N\max\{N-2,1\}$ where $\sigma_N$ is the measure of the unit Euclidean sphere, $\sigma_N\decl \frac{2\pi^{N/2}}{\Gamma(N/2)} $.
It follows that,
\begin{equation}\label{def:greenD}
\Gamma^x(y)=\begin{cases} \frac{C_N}{|x-y|^{N-2}}, &for\  N\ge 3,
\\-C_2\ln |x-y|,   &for\  N=2,
\end{cases} \quad for\ \ x\neq y.
\end{equation}

Let $x=(x',x_N)\in \HS$. We set $\hat x=(x',-x_N)$.
For $y\in \Rset^N_+\setminus\{ x\}$ we put,
\begin{equation}\label{def:G}
  G(x,y)\decl G^x(y) \decl \Gamma^x(y)-\Gamma^{\hat x}(y),
\end{equation}
and for $x\in\HS$ and $y'\in \DHS,$ we define,
\begin{equation}\label{def:K}
  K^x(y'):= \frac{C'_N x_N}{(x_N^2+|x'-y'|^2)^{N/2}}=
  \frac{\partial G^x}{\partial y_N}(y',0). %\quad where\ \   C_N'=\frac 2{\sigma_N}.
\end{equation}

By computation, we notice the following properties
of the Green's function $G$:
\begin{equation}
G^x(y)=G^y(x),\quad \forall x\neq y,
\end{equation}
\begin{equation}
G^x(y',0)=0,\quad \forall (x',x_N) \neq (y',0),
\end{equation}
\begin{equation}
 \frac{\partial G}{\partial y_N}(x',0,y)=0.
\end{equation}

For further applications, we need  a regularized version of the kernels $G$ and $K.$
For  $0\le \epsilon<1$, $x,y\in \CHS$, we define the 
regularized fundamental solution $\Gamma^x_\epsilon(y)$ 
 and the  regularized Green's function
$G_\epsilon^x(y),$ as follows
$$\Gamma^x_\epsilon(y)\decl\begin{cases}
 \frac{C_N}{(\epsilon^2+|x-y|^2)^{(N-2)/2}},&  N\ge 3,\\ \\
-\frac{C_2} 2\ln(\epsilon^2+|x-y|^2),  & N=2,
\end{cases}\quad
$$
and\ \ 
\begin{equation}\label{Geps}
G_\epsilon^x(y)\decl\Gamma^x_\epsilon(y) - \Gamma^{\hat x}_\epsilon(y).
\end{equation}
That is,
$$G_\epsilon^x(y)\decl\frac{C_N}{(\epsilon^2+|x-y|^2)^{(N-2)/2}}-\frac{C_N}{(\epsilon^2+|\hat x-y|^2)^{(N-2)/2}}
%= \frac{C_N}{\xp^{(N-2)/2}}-\frac{C_N}{\xg^{(N-2)/2}}
, \quad for\ N\ge 3,$$
$$G_\epsilon^x(y)\decl-\frac {C_2} 2\ln(\epsilon^2+|x-y|^2)+\frac {C_2} 2\ln(\epsilon^2+|\hat x-y|^2)
%=\frac{{C_2}}{2}\ln \frac{\xg}{\xp}
, \quad for\ N=2.$$
For $x\in\CHS$ and $y'\in\DHS$ we define the regularized Poisson kernel as
$$ K^x_\epsilon(y')\decl\frac{\partial G^x_\epsilon}{\partial y_N}(y',0)
=\frac{C'_N x_N}{(\epsilon^2+x_N^2+|x'-y'|^2)^{N/2}}.$$
These regularized kernels $G^x_\epsilon(y)$ and $K^x_\epsilon(y)$  help in dealing with potential singularities and are useful in various analytical settings.
Observe that, 
$$G_\epsilon^x(y)\nearrow G^x(y), \quad and \quad
K^x_\epsilon(y')\nearrow K^x(y')\ \ as\ \ \epsilon\to 0.$$

Let us introduce the auxiliary functions,
$$\xg\decl \epsilon^2+|\hat x-y|^2,\quad \xp\decl \epsilon^2+|x-y|^2.$$
It follows that  $G_\epsilon$ can be written as,
$$G_\epsilon^x(y)=
\frac{C_N}{\xp^{(N-2)/2}}-\frac{C_N}{\xg^{(N-2)/2}}, \quad for\ N\ge 3,$$
$$G_\epsilon^x(y)=\frac{{C_2}}{2}\ln \frac{\xg}{\xp}, \quad for\ N=2.$$
Notice that
\begin{equation}\label{a1-a2}
\xg-\xp=4x_N y_N.
\end{equation}
Therefore, in the case $N=2$, we have
\begin{equation}\label{Galt2}
G^x_\epsilon(y)=\frac{C_2}{2}\ln\frac{\xg}{\xp}=\frac{C_2}{2}\ln(1+\frac{\xg-\xp}{\xp})= \frac{C_2}{2}\ln(1+\frac{4x_Ny_N}{\xp}),
\end{equation}
while for $N\ge 3$, we have
\begin{eqnarray}
G^x_\epsilon(y)&=&\frac{C_N}{\xp^{(N-2)/2}}-\frac{C_N}{\xg^{(N-2)/2}}=
C_N\frac{\xg^{(N-2)/2}-\xp^{(N-2)/2}}{\xg^{(N-2)/2}\xp^{(N-2)/2}}\\&=&
C_N\frac{\xg^{N-2}-\xp^{N-2}}{\xg^{(N-2)/2}\xp^{(N-2)/2}(\xg^{(N-2)/2}
+\xp^{(N-2)/2})}\\
&=&C_N\frac{4x_Ny_N(\xg^{N-3}+\xg^{N-4}\xp+\dots+\xp^{N-3})}{\xg^{(N-2)/2}\xp^{(N-2)/2}(\xg^{(N-2)/2}
+\xp^{(N-2)/2})}.\label{Galt3}
\end{eqnarray}

For future reference, let us compute the derivatives of $G^x_\epsilon$. We have,
\begin{equation}\label{DNG}
\frac{\partial}{\partial y_N} G^x_\epsilon(y)=
\frac{C'_N}{2} \left(\frac{y_N+x_N}{\xg^{N/2}}-\frac{y_N-x_N}{\xp^{N/2}}\right),
\end{equation}
and for $j=1,\dots,N-1$,
$$\frac{\partial}{\partial y_j} G^x_\epsilon(y)=
\frac{C'_N}{2} \left(\frac{1}{\xg^{N/2}}-\frac{1}{\xp^{N/2}}\right)(y_j-x_j).$$
So that,
\begin{eqnarray}
  |\nabla G^x_\epsilon|^2&=&\left(\frac{C'_N}{2}\right)^2
      \abs{\frac{y-\hat x}{\xg^{N/2}} -\frac{y-x}{\xp^{N/2}}}^2\\
      &=&     \left(\frac{C'_N}{2}\right)^2\left[
     \abs{\frac{1}{\xg^{N/2}} -\frac{1}{\xp^{N/2}} }^2  |x'-y'|^2+
     \left(\frac{y_N+x_N}{\xg^{N/2}}-\frac{y_N-x_N}{\xp^{N/2}}\right)^2\right].\label{eq:GRG}
\end{eqnarray}

This detailed computation of the gradients and their magnitudes provides a clear understanding of the behavior of the regularized Green's function  and its derivatives. This will be useful in further analytical work involving these functions.

\begin{proposition}\label{prop1} Let $0\le\epsilon<1$. 
\item[1.] For $y,x\in\HS$, $x\neq y$, we have
\begin{equation}\label{estGabove}
0\le G^x_\epsilon(y) \le C'_N \frac{x_Ny_N}{(\epsilon^2+|x-y|^2)^{N/2}}.\end{equation}

\item[2.] We have,
\begin{equation}\label{estGasym}
\lim_{|y|\to\infty} G^x_\epsilon(y) \frac{|y|^N}{y_N}=C'_Nx_N,\end{equation}
uniformly with respect to $\epsilon$ and $x$ in bounded sets.

\item[3.] We have,
\begin{equation}\label{estKasym}
\lim_{|y'|\to\infty} K^x_\epsilon(y') {|y'|^N}=C'_Nx_N,\end{equation}
uniformly with respect to $\epsilon$ and $x$ in bounded sets.
\item[4.] For $y,x\in\HS$, $x\neq y$, we have
\begin{eqnarray}
 |\nabla G_\epsilon^x|^2 (y) &\le&\left(\frac{C'_N}{2}\right)^2\frac{x_N^2}{\xp^N}
\left[ \left(1+(\frac{\xp}{\xg})^{N/2}+2N\frac{y_N^2}{\xg}\right)^2+4N^2\frac{y_N^2 |x'-y'|^2}{\xg^2} \right]\label{est:GRabove1}\\
&\le& {C'_N}^2\frac{x_N^2}{(\epsilon^2+|x-y|^2)^{N}}\left[ (1+N)^2+N^2\right],\label{est:GRabove2}\\
 |\nabla G_\epsilon^x|^2(y)  &\ge &\left(\frac{C'_N}{2}\right)^2\frac{x_N^2}{(\epsilon^2+|x-y|^2)^{N}}.\label{est:GRbelow}
\end{eqnarray}

\item[5.] 
 For $y,x\in\CHS$, $x\neq y$, we have
\begin{equation}\label{estDNGabove}
|\partial_N G^x_\epsilon(y)| \le \frac{C'_N}{2}\frac{x_N}{\xp^{N/2}} 
  \left[1+\left(\frac{\xp}{\xg}\right)^{N/2} + 2N\frac{y_N^2}{\xg}\right]\le 
   (N+1)C'_N\frac{x_N}{(\epsilon^2+|x-y|^2)^{N/2}}.
\end{equation}
\end{proposition}

\begin{proof} Proof of claim 1. We begin by observing that since $\xg\ge\xp$, we have $ G^x_\epsilon\ge 0$.

In the case $N=2$, by the concavity of the logarithmic function, we obtain:
$$G^x_\epsilon=\frac{C_2}{2}\ln\frac{\xg}{\xp}=\frac{C_2}{2}\ln(1+\frac{\xg-\xp}{\xp})\le \frac{C_2}{2}\frac{\xg-\xp}{\xp} = \frac{C_2}{2}\frac{4x_Ny_N}{\xp},$$
which proves \eqref{estGabove}.

Let $N\ge 3$. By the convexity of the function, $t\mapsto t^{-\alpha}$ for $t>0$ and $\alpha>0$, we have:
\begin{equation}\label{conv1}
  \frac{1}{s^\alpha}-\frac{1}{t^\alpha}\le\alpha\frac{1}{s^{\alpha+1}}(t-s),\qquad for\ \ t,s>0.
\end{equation}
Applying \eqref{conv1} with $\alpha=(N-2)/2$, $s=\xp$ and $t=\xg$, we deduce:
$$  G^x_\epsilon=C_N\left(\frac{1}{\xp^{(N-2)/2}}-\frac{1}{\xg^{(N-2)/2}}\right)\le C_N \frac{N-2}{2}\frac{1}{\xp^{N/2}}(\xg-\xp)=
C_N(N-2)2\frac{x_Ny_N}{\xp^{N/2}}.
$$
% Since
%$\xg\ge \xp$,  from \eqref{Galt3} we deduce
%$$\begin{aligned}
%0\le G^x_\epsilon&=
%C_N\frac{4x_Ny_N(\xg^{N-3}+\xg^{N-4}\xp+\dots+\xp^{N-3})}{\xg^{(N-2)/2}%\xp^{(N-2)/2}(\xg^{(N-2)/2}
%+\xp^{(N-2)/2})}\\
%&\le C_N\frac{4x_Ny_N (N-2) \xg^{N-3}}{(\xg\xp)^{(N-2)/2} 2 \xp^{(N-2)/2}}
% =C'_N\frac{x_Ny_N }{\xp^{N/2}}\frac{ \xg^{N-3}} {(\xg\xp)^{(N-2)/2}  \xp^{-1}},
%\end{aligned}
%$$
This completes the proof of \eqref{estGabove} for $N\ge 3.$

\noindent Proof of claim 2. Let $C\subset\HS$ be bounded and let $x\in C$. For
$|y|>2 \sup \{|x|:x\in C	\}$, we have,
$$\sqrt{1+|y-\hat x|^2}\ge \sqrt\xg\ge\sqrt \xp\ge |y-x|\ge |y|-|x|\ge |y|/2.
$$

Next, let $N=2$. 
Since $\frac{4x_Ny_N}{\xp}\to 0$ uniformly with respect to $x\in C$ and $\epsilon,$
whenever $|y|\to\infty$, and taking into account that
$\ln(1+t)\cong t$ as $t\to 0$,\footnote{$a(t)\cong b(t)$ 
as $t\to t_0$ means $\lim_{t\to t_0} a(t)/b(t)=1$.}  
from \eqref{Galt2} it follows that
$$G^x_\epsilon(y)\cong \frac{C_2}{2} \frac{4x_Ny_N}{\xp}\cong 2C_2 x_N\frac{y_N}{|y|^2},\quad as \ \ |y|\to\infty,
$$
completing the proof of the claim.

Next we consider the case $N\ge 3$.
From \eqref{Galt3} we have,
$$\begin{aligned}
G^x_\epsilon(y)\frac{|y|^N}{y_N}=& 4C_Nx_N\frac{(\xg^{N-3}+\xg^{N-4}\xp+\dots+\xp^{N-3})}{|y|^{2N-6}}\frac{|y|^{2N-4}}{(\xg\xp)^{(N-2)/2}}\frac{|y|^{N-2}}{\xg^{(N-2)/2}
+\xp^{(N-2)/2}}\\
\cong&  4C_Nx_N (N-2)\frac12,\qquad as\ \  |y|\to\infty,
\end{aligned}
$$
this completes the proof of claim 2.

\noindent  Proof of claim 3. Proceeding as in the proof of claim 2., we have 
$$K^x_\epsilon(y')|y'|^N=C'_Nx_N\frac{|y'|^N}{(\epsilon^2+x_N^2+|x'-y'|^2)^{N/2}}\cong C'_Nx_N,\quad as \ \ |y'|\to\infty.$$

\noindent Proof of claim 4. and claim 5.
From the expression of the derivative of $G^x_\epsilon,$ it follows that we must estimate the quantity,
$\xp^{-N/2}-\xg^{-N/2}$. 
We observe that since the function $t^\alpha$ is convex for $\alpha\ge 1$, we have
\begin{equation}\label{conv}
t^\alpha-s^\alpha\le\alpha (t-s) t^{\alpha -1},\qquad for\ \ t,s>0.
\end{equation}
Taking into account \eqref{a1-a2}, 
the latter inequality, with $t=\xg$, $s=\xp$ and $\alpha=N/2$, produces
$$\frac{1}{\xp^{N/2}}-\frac{1}{\xg^{N/2}}=\frac{\xg^{N/2}-\xp^{N/2}}{\xg^{N/2}\xp^{N/2}}\le 
\frac N2(\xg-\xp)\frac{\xg^{N/2-1}}{\xg^{N/2}\xp^{N/2}}=
2N\frac{x_Ny_N}{\xg\xp^{N/2}}.
$$

Next from \eqref{DNG} we have,
$$\begin{aligned}
|\partial_N G^x_\epsilon(y)|=&
\frac{C'_N}{2} \left(y_N(\frac{1}{\xp^{N/2}}-\frac{1}{\xg^{N/2}})+x_N (\frac{1}{\xp^{N/2}}+\frac{1}{\xg^{N/2}})\right)\\
\le&  \frac{C'_N}{2} \left( y_N2N\frac{x_Ny_N}{\xg\xp^{N/2}} + x_N (\frac{1}{\xp^{N/2}}+\frac{1}{\xg^{N/2}})\right).
\end{aligned}
$$
This proves the first inequality in \eqref{estDNGabove},
while the second one is due to the fact that,
$\xg\ge \xp$ and $\xg\ge y_N^2$. This concludes the proof of claim 5.

Similarly form \eqref{eq:GRG} we have,
$$\begin{aligned}
\frac{ |\nabla G^x_\epsilon|^2}{({C'_N}/{2})^2}=&
     \abs{\frac{1}{\xg^{N/2}} -\frac{1}{\xp^{N/2}} }^2  |x'-y'|^2+
     \left(y_N(\frac{1}{\xp^{N/2}}-\frac{1}{\xg^{N/2}})+x_N (\frac{1}{\xp^{N/2}}+\frac{1}{\xg^{N/2}})
     \right)^2\\
     \le& 4N^2\frac{x_N^2y_N^2}{\xg^2\xp^{N}}|x'-y'|^2+
     \frac{x_N^2}{\xp^{N}} 
  \left(1+\left(\frac{\xp}{\xg}\right)^{N/2} + 2N\frac{y_N^2}{\xg}\right)^2,
\end{aligned}
$$
which is \eqref{est:GRabove1}. 
The remaining inequality \eqref{est:GRabove2} is a consequence of the fact that,
$\xg\ge \xp$, $\xg\ge y_N^2$ and $\xg\ge  |x'-y'|^2$. 

Second, we prove the estimate \eqref{est:GRbelow}. Indeed we have,
$$\begin{aligned}
\frac{ |\nabla G^x_\epsilon|^2}{({C'_N}/{2})^2}=&
     \abs{\frac{1}{\xg^{N/2}} -\frac{1}{\xp^{N/2}} }^2  |x'-y'|^2+
     \left(y_N(\frac{1}{\xp^{N/2}}-\frac{1}{\xg^{N/2}})+x_N (\frac{1}{\xp^{N/2}}+\frac{1}{\xg^{N/2}})
     \right)^2\\
     \ge& \left(x_N (\frac{1}{\xp^{N/2}}+\frac{1}{\xg^{N/2}})
     \right)^2\ge x_N^2\frac{1}{\xp^{N}}.
  \end{aligned}
$$   
\end{proof}

\section{Examples of functions satisfying {\rhz}}\label{sec:R0}
 Here we summarize examples of functions satisfying \rhz\!\!.
\begin{remark}\label{rem:ring} We notice that a function $u\in L^1_{loc}(\CHS)$ can be extended to the whole space  $\mathbb{R}^N$
by setting $u(x)=0$ for $x\notin \CHS$. We will still  denote this extension with $u$. If $u$
satisfies the ring condition 
\begin{equation}\label{ringentireSQ}
\liminf\frac{1}{R^N}\int_{R<\abss{x-y}<2R} |u(y)|dy=0, \qquad for\  a.e.\ x\in\HS,
\end{equation}
 then $u$ satisfies also \rhz\!\!.
 It is easy to see that \eqref{ringentireSQ} is satisfied by any function belonging to the following global spaces:
 $L^p(\HS)$, the Marcinkiewicz spaces $L^p_w(\HS)$, and to the Campanato-Morrey spaces $M^p_q(\HS)$ (see \cite{CDM}).
 This implies that such a function fulfills \hbox{\rhz\!\!\!\!~.} 
In effect, essentially bounded functions fill  \rhz but not \eqref{ringentireSQ}.
\end{remark}

\begin{theorem} \label{teo:ringLP}
1. Let $1\le p\le\infty$. 
If $u\in L^p(\HS)$, then   $u$ satisfies \rhz\!\!.

2. Let $1\le p<\infty$.  If $x_N|u|^p\in L^1(\HS)$, then   $u$ satisfies \rhz\!\!.

3. If $x_N|u|\in L^\infty(\HS)$, then   $u$ satisfies \rhz\!\!.

4. If $u\in L^1_{loc}(\CHS)$ and $|u(x)|\le c|x|^p$ for a suitable $c>0$, $0<p<1$ and 
$|x|$ large, then   $u$ satisfies \rhz\!\!. 

5.  If $u\in M^p_q(\HS)$ with $1\le q\le p$, then   $u$ satisfies \rhz\!\!.
\end{theorem}

\begin{proof}We begin by proving  statements 2 and 4, since they imply the others.
Let $x\in\HS$. For $R>0$, define
$A^*_R(x)\decl B^*_{2R}(x)\setminus B^*_R(x).$
See Figure \ref{fig1} at page \pageref{fig1}.

\noindent Proof of claim 2. From our assumption, we know that
$$ \lim_{R\to \infty}\int_{A^*_R(x)}y_N|u(x)|^p dx=0.$$

This directly implies the claim for $p=1$.
Now consider the case $p>1$. By H\"older inequality with exponent $p$, we have
$$\begin{aligned}
 \frac{1}{R^{N+2}}\int_{A^*_R(x)}y_N|u(y)| dy
=&\frac{1}{R^{N+1}}\int_{A^*_R(x)}|u(y)|\frac{y_N}{R} dy\\
\le& \frac{1}{R^{N+1}}
\left(\int_{A^*_R(x)}|u(y)|^p\frac{y_N}{R} dy
\right)^{1/p}
\left(\int_{A^*_R(x)}\frac{y_N}{R} dy
\right)^{1/p'}\\
\le &
 \frac{R^{N/p'}}{R^{N+1}}
\left(\int_{A^*_R(x)}|u(y)|^p\frac{y_N}{R} dy
\right)^{1/p}\to 0,
\end{aligned}
$$ 
and the claim follows.

\noindent Proof of claim 4. Let $R_0$ be such that $|u(x)|\le c |x|^p\le c(|x'|^p+|x_N|^p)$ for $|x|>R_0$.
We have
$$\begin{aligned}
\int_{B^*_R\setminus B^*_{R_0}} y_N|u(y)| dy\le 
c \int_{B^*_R\setminus B^*_{R_0}} y_N (|y'|^p+|y_N|^p) dy\le 
c  \int_{B^*_R} y_N (|y'|^p+|x_N|^p) dy \\
\le \int_{|y'|<R} dy' \int_0^R  y_N (|y'|^p+|y_N|^p) dy_N
\le c \int_{|y'|<R} ( |y'| R^2 + R^{p+2}) dy' 
\le c R^{N+1+p}.
\end{aligned}
$$
Therefore, we can complete the proof by noticing that
$$ \frac{1}{R^{N+2}} \int_{B^*_R\setminus B^*_{R_0}} y_N|u(y)| dy\le
 c R^{p-1}\to 0,\ \ as\ \ R\to\infty.$$

\noindent Proof of claims  1. and 5. The case $1\le p<\infty$ is contained in Remark \ref{rem:ring}.
While the case $p=\infty$ follows from point 4.
%$$ \frac{1}{R^{N+2}}\int_{A^*_R(y)}y_N|u(x)| dy\le 
% \frac{||u||_\infty}{R^{N+1}}\int_{A^*_R(y)} dy\le c  \frac{||u||_\infty}{R}\to 0.
%$$ 

\noindent Proof of claim 3. It follows by using the same argument as above.
\end{proof}

%\medskip

Another example of functions  which naturally meet the condition \rhz is the following.

\begin{theorem} \label{teo:power} Let $1< q<\infty$. Let $u\in L^q_{loc}(\HS)$ be a weak solution
%\footnote{The formal definition of weak solution will be introduced in Section \ref{sec3} below.} 
of
\begin{equation}\label{dispower}
 \left\{\begin{array}{ll}
-\Delta u \ge |u|^q, & \mathrm{on }\ \HS, \cr
u\ge 0, & \mathrm{on }\ \DHS. \cr
\end{array}
 \right.
\end{equation}
Then  $u$ satisfies \rhz\!\!.
\end{theorem}
For a proof of this theorem see \cite{damiII}.

\end{appendices}

\medskip 

\noindent \textbf{Data Availability.} This manuscript has no associated data

\bibliographystyle{amsplain}

\begin{thebibliography}{10}


\bibitem{AG} {\sc D.H. Armitage} and {\sc S.J. Gardiner}, 
Classical potential theory. London: Springer, 2001.

\bibitem{axl} {\sc S. Axler, P. Bourdon} and {\sc W. Ramey},
Harmonic function theory. 2nd ed. New York, NY: Springer, 2001.



\bibitem{BNV} {\sc H. Berestycki, L. Nirenberg} and 
{\sc S.R.S. Varadhan}, The principal eigenvalue and maximum principle for second-order elliptic operators in general domains.
 \textit{Comm. Pure Appl. Math.}
\textbf{47} (1994), 47--92. %doi:10.1002/cpa.3160470105 


\bibitem{BLU}{\sc A. Bonfiglioli, E. Lanconelli} and {\sc  F. Uguzzoni}, Stratified Lie groups and potential theory for their sub-Laplacians. Springer Monographs in Mathematics. Springer, Berlin, 2007. xxvi+800 pp. 


\bibitem{bou} {\sc N. Bourbaki},
Elements of mathematics. Integration I: Chapters 1--6. Translated from the 1959, 1965 and 1967 French originals by Sterling K. Berberian. Berlin: Springer, 2004.


 
\bibitem{CDM} {\sc G. Caristi, L. D'Ambrosio} and {\sc E.  Mitidieri}, 
Representation formulae for solutions to some classes of higher order systems and related Liouville theorems,
 \textit{Milan J. Math. }\textbf{76} (2008), 27--67.
 
\bibitem{DG22} {\sc L. D'Ambrosio} and {\sc M. Ghergu},
{Representation formulae for nonhomogeneous differential operators and applications to PDEs},
 \textit{J. Differential Equations} {\bf 317} (2022), 706--753.



\bibitem{DM19}{\sc L.~D'Ambrosio} and {\sc E.~Mitidieri},
{Representation formulae of solutions of second order elliptic
inequalities},
\textit{Nonlinear Anal.} {\bf 178} (2019), 310--336.

\bibitem{damiII}{\sc L.~D'Ambrosio} and {\sc E.~Mitidieri}, \newblock Liouville theorems of semilinear elliptic inequalities in a half-space, (2025), 1--29, preprint.

\bibitem{DMP} {\sc L. D'Ambrosio, E. Mitidieri} and {\sc S. I.  Pohozaev}, {Representation formulae and inequalities for solutions of a class of second order partial differential equations},
\textit{Trans. Amer. Math. Soc.} \textbf{358} (2006), 893--910.



\bibitem{dau-lio80}{\sc R. Dautray} and {\sc J.-L. Lions},
Mathematical analysis and numerical methods for science and technology. Volume 1: Physical origins and classical methods. With the collaboration of Ph. B\'enilan, M. Cessenat, A. Gervat, A. Kavenoky, H. Lanchon. Transl. from the French by I.N. Sneddon. Berlin: Springer-Verlag, 1990.

\bibitem{evagar} {\sc L.C. Evans} and {\sc R.F. Gariepy},  Measure theory and fine properties of functions. Boca Raton: CRC Press, 1992.



\bibitem{giu84}{\sc E. Giusti},
 Minimal surfaces and functions of bounded variation. 
Monographs in Mathematics, Vol. 80. Boston-Basel-Stuttgart: Birkhäuser, 1984. 


\bibitem{H56}  {\sc A. Huber}, 
{On functions subharmonic in a half-space},
 \textit{Trans. Amer. Math. Soc.} {\bf 82} (1956), 147-159.  

\bibitem{leoni} {\sc G. Leoni}, A first course in Sobolev spaces. Providence, RI: American Mathematical Society (AMS),  2009.

\bibitem{marver} {\sc M. Marcus} and {\sc L. V\'eron}, 
 Nonlinear second order elliptic equations involving measures. Berlin: de Gruyter, 2014.
 
\bibitem{zha16} {\sc Y.H. Zhang, G.T. Deng} and {\sc T. Qian},
Integral representations of a class of harmonic functions in the half space,
\textit{J. Differential Equations}
 \textbf{260} (2016), 923--936.


 
\end{thebibliography}

\end{document}